\documentclass[10pt]{amsart}
\usepackage{amscd,amsxtra,amssymb,mathrsfs, bbm}
\usepackage{stmaryrd}

\newtheorem{theorem}{Theorem}[section]
\newtheorem{corollary}[theorem]{Corollary}
\newtheorem{lemma}[theorem]{Lemma}
\newtheorem{proposition}[theorem]{Proposition}

\theoremstyle{definition}
\newtheorem{definition}[theorem]{Definition}
\newtheorem{remark}[theorem]{Remark}

\DeclareMathAlphabet{\mathpzc}{OT1}{pzc}{m}{it}

\DeclareMathOperator{\rk}{rk}

\newcommand{\ev}{\mathsf{ev}}
\renewcommand{\Im}{\mathsf{Im}}
\renewcommand{\dim}{\mathsf{dim}}

\DeclareMathOperator{\Coh}{\mathsf{Coh}}

\DeclareMathOperator{\Pic}{Pic}

\DeclareMathOperator{\Rep}{\mathsf{Rep}}

\DeclareMathOperator{\Hom}{\mathsf{Hom}}

\DeclareMathOperator{\Tor}{\mathsf{Tor}}

\DeclareMathOperator{\Ext}{\mathsf{Ext}}

\DeclareMathOperator{\GL}{\mathsf{GL}}
\DeclareMathOperator{\Aut}{\mathsf{Aut}}

\DeclareMathOperator{\End}{\mathsf{End}}

\input xy
\xyoption{all}

\newcommand{\kk}{k}

\newcommand{\FF}{\mathbb{F}}
\newcommand{\GG}{\mathbb{G}}
\renewcommand{\mod}{\mathsf{mod}}

\renewcommand{\AA}{\mathbb{A}}
\newcommand{\DD}{\mathbb{D}}
\renewcommand{\SS}{\mathbb{S}}

\newcommand{\PP}{\mathbb{P}}

\newcommand{\RR}{\mathbb{R}}
\newcommand{\LL}{\mathbb{L}}
\newcommand{\HH}{\mathbb{H}}

\newcommand{\QQ}{\mathbb Q}
\newcommand{\CC}{\mathbb C}
\newcommand{\ZZ}{\mathbb Z}

\newcommand{\kF}{\mathcal{F}}

\newcommand{\kO}{\mathcal{O}}

\newcommand{\kP}{\mathcal{P}}
\newcommand{\kQ}{\mathcal{Q}}

\newcommand{\kT}{\mathcal{T}}

\newcommand{\lar}{\longrightarrow}
\newcommand{\overr}{\overrightarrow}

\newcommand{\cA}{\mathsf{A}}
\newcommand{\cB}{\mathsf{B}}
\newcommand{\cI}{\mathsf{J}}
\newcommand{\cR}{\mathsf{R}}

\begin{document}

\title[Two descriptions of $U_v(\widehat{\mathfrak{sl}}_2)$]{Two descriptions of  the quantum
 affine algebra $U_v(\widehat{\mathfrak{sl}}_2)$ via Hall algebra approach}

\author{Igor Burban}
\address{
Mathematisches Institut,
Universit\"at Bonn,
Endenicher Allee 60,
D-53115 Bonn,
Germany
}
\email{burban@math.uni-bonn.de}

\author{Olivier Schiffmann}
\address{
D\'epartement de Math\'ematiques \\
Universit\'e  Paris Sud \\
B\^atiment 425 \\
91405 Orsay Cedex \\
 France
}
\email{Olivier.Schiffmann@math.u-psud.fr}

\subjclass[2000]{Primary 16T20, 18E30,  Secondary 17B37}
\keywords{Quantum affine algebras, Hall algebras, triangulated and derived categories, tilting}

\begin{abstract}
We compare the reduced Drinfeld doubles of the composition subalgebras of
the category of representations of the Kronecker quiver $\overr{Q}$ and of the category of
coherent sheaves on $\PP^1$. Using this approach,
we show that the   Drinfeld--Beck isomorphism for the quantized enveloping algebra
$U_v(\widehat{\mathfrak{sl}}_2)$ is a  corollary of an equivalence between the
derived categories
$D^b\bigl(\Rep(\overr{Q})\bigr)$ and $D^b\bigl(\Coh(\PP^1)\bigr)$. This technique allows
to reprove several  results on the integral form of $U_v(\widehat{\mathfrak{sl}}_2)$.
\end{abstract}

\maketitle

\vspace{.1in}

\section{Introduction}

\noindent
In this article, we study the  relation between the composition algebra of the
category of representations of the Kronecker quiver
$$
\overr{Q}
=
\xymatrix{
\bullet  \ar@/^/[rr] \ar@/_/[rr] & & \bullet
}
$$
and the composition algebra
of the category
of coherent sheaves on the  projective line $\PP^1$.
As it was shown by Ringel \cite{Ringel} and Green \cite{Green}, the generic composition
algebra of the category of representations of $\overr{Q}$ is isomorphic to the positive
part of the quantum  affine algebra  $U_v(\widehat{\mathfrak{sl}}_2)$ written
in the terms of Drinfeld--Jimbo generators.

On the other side, as  it was discovered by Kapranov  \cite{Ka} and extended  by Baumann and Kassel
 \cite{BK}, the Hall algebra of the category of coherent sheaves on a
projective line $\PP^1$ is closely related
with Drinfeld's new realization
$U_v\bigl(\mathfrak{L sl}_2\bigr)$ of the quantized enveloping algebra  of $\widehat{\mathfrak{sl}}_2$
\cite{Drinfeld}.
Since then, this subject drew attention of many authors, see for example
\cite{Zhang, Szanto, McGerty,XiaoZhang, SchiffmannInvent, OlivierNotes}.

In this article, we work out this important observation a step further and show that
the Drinfeld--Beck isomorphism $U_v\bigl(\widehat{\mathfrak{sl}}_2\bigr) \rightarrow
U_v\bigl(\mathfrak{L sl}_2\bigr)$ (see \cite{Drinfeld, Beck, DingFrenkel, Jing})
can be viewed as a corollary of the derived equivalence
$D^b(\Rep(\overr{Q})\bigr) \rightarrow D^b\bigl(\Coh(\PP^1)\bigr)$. The understanding of
this isomorphism is of great  importance for the representation theory of $U_v\bigl(\widehat{\mathfrak{sl}}_2\bigr)$ and its applications in mathematical physics, see for example
\cite{JimboMiwa}. Indeed, since $U_v\bigl(\widehat{\mathfrak{sl}}_2\bigr)$
is a Hopf algebra, the category of its finite-dimensional representations has a structure of a
tensor category. However, in order to describe such representations themselves, it is frequently
more convenient to work with Drinfeld's new realization  $U_v\bigl(\mathfrak{L sl}_2\bigr)$.
Using the fact that an equivalence of triangulated categories commutes with Serre functors, we also  show  that the Drinfeld-Beck isomorphism $U_v\bigl(\widehat{\mathfrak{sl}}_2\bigr) \rightarrow
U_v\bigl(\mathfrak{L sl}_2\bigr)$
commutes with the Coxeter transformation acting on both sides. Finally, applying
 the technique of stability conditions,
we  reprove several known technical statements on the integral form of $U_v\bigl(\mathfrak{L sl}_2\bigr)$.

\vspace{1mm}
\noindent
\emph{Notation}. Throughout the paper, $\kk = \FF_q$ is a finite field with $q$ elements
and   $\widetilde\QQ =  \widetilde{\QQ}_q = \QQ[v, v^{-1}]/(v^{-2}- q) \cong \QQ[\sqrt{q}]$.
Next, $\kP$ is the set of all integers of the form  $p^t$, where $p$ is a prime number
and $t \in \ZZ_+$.  For a positive integer $n$ we set $[n] = [n]_v =
\frac{\displaystyle v^n - v^{-n}}{\displaystyle v - v^{-1}}$
and $[n]! = [1] \dots  [n]$. We denote by $R$ the subring
of the field $\QQ(v)$ consisting of the rational functions having poles only at $0$
or at roots of $1$.  For the  affine  Lie algebra $\mathfrak{g} = \widehat{\mathfrak{sl}}_2$
we denote by $U_v(\mathfrak{g})$ its  quantized enveloping algebra over the ring $R$, whereas
$U_q(\mathfrak{g}) = U_v(\mathfrak{g}) \otimes_{R} \widetilde{\QQ}$.

\vspace{2mm}
\noindent
\emph{Acknowledgement}.   The research of the first-named
author was supported by the DFG project Bu--1866/2--1. Some parts of this work were done
during the authors stay at  the Mathematical Research Institute in Oberwolfach within
the ``Research in Pairs'' programme.

\section{Hereditary categories, their Hall algebras and Drinfeld doubles}

Let $\cA$ be an
essentially small hereditary abelian $\kk$-linear category such that for all objects $M, N \in \mathrm{Ob}(\cA)$
the $\kk$--vector spaces $\Hom_\cA(M, N)$ and $\Ext^1_\cA(M, N)$ are finite dimensional. In what follows,
we shall call such a category \emph{finitary}.
Let $\cI = {\cI}_\cA := \bigl(\mathrm{Ob}(\cA)/\cong\bigr)$ be the set of isomorphy classes of objects in $\cA$.
For an object $X \in \mathrm{Ob}(\cA)$,  we denote by $[X]$ its image in $\cI$. Fix the following notations.

\begin{itemize}
\item For any object $X \in \mathrm{Ob}(\cA)$ we set $a_X = \bigl|\Aut_\cA(X)\bigr|$.
\item For any three objects $X, Y, Z \in \mathrm{Ob}(\cA)$ we denote
$$P^{Z}_{X, Y} = \left| \left\{(f, g) \in \Hom_{\cA}(Y, Z) \times \Hom_{\cA}(Z, X) \, \, \left| \, \,
0 \to Y \stackrel{f}\rightarrow  Z \stackrel{g}\rightarrow X \to 0 \right.
\mbox{\, \, is exact} \right\}\right|.
$$
\item Finally, we put  $F^{Z}_{X, Y} = \frac{\displaystyle P^{Z}_{X, Y}}{\displaystyle a_X \cdot a_Y}$.
\end{itemize}
Note that the numbers $a_Z$, $P^{Z}_{X, Y}$,  $F^{Z}_{X, Y}$ and $\frac{\displaystyle P^{Z}_{X, Y}}{\displaystyle a_Z}$ depend only on the isomorphy classes of $X, Y, Z$ and are integers.

Let $K = K_0(\cA)$ be the K--group of $\cA$. For an object $X \in \mathrm{Ob}(\cA)$,  we denote by
$\bar{X}$ its image in $K$. Next, let $\langle \,-\,, \,-\,\rangle: K \times K \rightarrow \ZZ$ be the
Euler form:
$$
\langle \bar{X}, \bar{Y}\rangle = \dim_{\kk} \Hom_\cA(X, Y) - \dim_\kk \Ext^1_\cA(X, Y)
$$
and $(\,-\,, \,-\,): K \times K \rightarrow \ZZ$ its symmetrization:
$$
(\alpha, \beta) = \langle\alpha, \beta\rangle + \langle\beta, \alpha\rangle, \quad \alpha, \beta
\in K.
$$
Following  Ringel \cite{Ringel},  one  can attach to a finitary hereditary category $\cA$ an associative
algebra  $H(\cA)$ called the  extended twisted \emph{Hall algebra} of $\cA$,
defined over the field
$\widetilde{\QQ}$. As a vector space over $\widetilde{\QQ}$, we have
$$
\bar{H}(\cA) := \bigoplus\limits_{[Z] \in \cI} \widetilde{\QQ}[Z] \quad \mbox{and} \quad
H(\cA) :=  \bar{H}(\cA) \otimes_{\widetilde{\QQ}} \widetilde{\QQ}[K].
$$
For a class $\alpha \in K$ we denote by $K_\alpha$ the corresponding element in the group algebra $\widetilde{\QQ}[K]$. Then we have: $K_\alpha \circ  K_\beta = K_{\alpha + \beta}$.
Next, for   $[X], [Y] \in \cI$ the product $\circ$ is defined to be
$$
[X] \circ [Y] = \sqrt{\frac{|\mathrm{Hom}(X,Y)|}{|\Ext^1(X,Y)|}}
\sum\limits_{[Z] \in \cI} F^{Z}_{X, Y}\; [Z] = v^{-\langle \bar{X}, \bar{Y}\rangle} \sum\limits_{[Z] \in \cI} F^{Z}_{X, Y}\; [Z].
$$
Finally, for any $\alpha \in K$ and $[X] \in \cI$ we have
$$
K_\alpha \circ [X] = v^{-(\alpha, \bar{X})} [X] \circ K_\alpha.
$$
As it was shown in \cite{Ringel},  the product $\circ$ is associative and the element
$1 := [0] \otimes K_0$ is the unit element.  In what follows, we
shall use the notation $[X]K_\alpha$ for the  element $[X] \otimes K_\alpha  \in
H(\cA)$.

Let $\cA$ be a finitary \emph{finite length}  hereditary category over $\kk$. By a result of
Green \cite{Green}, the Hall algebra $H(\cA)$ has a natural bialgebra structure, see
also \cite{RingelGreen}.
The comultiplication $\Delta: H(\cA) \rightarrow H(\cA)
\otimes_{\widetilde{\QQ}} H(\cA)$
and the counit $\eta: H(\cA) \rightarrow \widetilde{\QQ}$  are
 given by the following formulae:
$$
\Delta\bigl([Z] K_\alpha\bigr) = \sum\limits_{[X], \, [Y] \in \cI}
v^{-\langle \bar{X}, \bar{Y}\rangle}
\frac{\displaystyle P^{Z}_{X, Y}}{a_Z} [X] K_{\bar{Y} + \alpha} \otimes [Y] K_\alpha \quad
\mbox{and} \quad
\eta\bigl([Z] K_\alpha\bigr) = \delta_{Z, 0}.
$$
Moreover, as it was shown by Xiao \cite{Xiao}, the Hall algebra $H(\cA)$  is also a Hopf algebra.
Finally, there is a  pairing $(\, -\,,\,-\,): H(\cA) \times H(\cA) \rightarrow
\widetilde{\QQ}$ introduced by Green \cite{Green},  given by the expression
$$
\Bigl([X] K_\alpha, [Y] K_\beta\Bigr) = v^{-(\alpha, \beta)}
\frac{\displaystyle \delta_{X, Y}}{\displaystyle a_X}.
$$
This pairing is non-degenerate on $\bar{H}(\cA)$ and symmetric. Next, it  satisfies the following properties:
$$
(a \circ b, c) = \bigl(a \otimes b, \Delta(c)\bigr) \quad \mbox{and} \quad  (a, 1) = \eta(a)
$$
for any $a,b,c \in H(\cA)$. In other words, it is a bialgebra  pairing.

\begin{remark}\label{R:TopolBialg}
If  $\cA$ is not a category of finite length (for instance, if it is the category of coherent sheaves
on a projective curve) then the Green's pairing
$(\,-\,,\,-\,): H(\cA) \times H(\cA) \rightarrow
\widetilde{\QQ}$ is still a well-defined symmetric bilinear pairing.
However, the comultiplication
$\Delta([X])$  is possibly an  infinite sum.
Nevertheless, it is possible to introduce a completed tensor product
$H(\cA) \widehat\otimes H(\cA)$ (which is a $\widetilde{\QQ}$--algebra) such
that $\Delta: H(\cA) \rightarrow H(\cA) \widehat\otimes H(\cA)$ is an algebra homomorphism
and $(\Delta \otimes \mathbbm{1})\circ  \Delta = (\mathbbm{1} \otimes \Delta) \circ \Delta$.
Moreover,
for any elements $a, b, c \in H(\cA)$
the expression $\bigl(a \otimes  b, \Delta(c)\bigr)$
 takes a  finite value  and the equalities
 $(a \circ b, \, c) = \bigl(a \otimes b, \Delta(c)\bigr), (a, 1) = \eta(a)
$ are fulfilled. In such a situation we say that $H(\cA)$ is a \emph{topological bialgebra},
see \cite[Appendix B]{HallEll} for further details.
\end{remark}

From now on, let $\cA$ be an arbitrary $\kk$--linear hereditary finitary abelian category.
Consider the root category $\cR(\cA) = D^b(\cA)/[2]$. Note that
$\cR(\cA)$ has a structure of  a triangulated category such that the canonical functor
$D^b(\cA) \rightarrow \cR(\cA)$ is exact, see \cite[Section 7]{PengXiao}. Moreover,
 any object of $\cR(\cA)$ splits into
a direct sum $X^+ \oplus X^-$, where $X^+ \in \mathrm{Ob}(\cA)$ and
$X^- \in \mathrm{Ob}(\cA)[1]$.

Our next goal is to introduce the reduced Drinfeld double of the topological bialgebra
$H(\cA)$. Roughly speaking (although, not completely correctly), it is
an analogue of the Hall algebra, attached to the triangulated category $\cR(\cA)$.
To define it,  consider the pair of   algebras $H^\pm(\cA)$, where  we use the notation
$$
H^+(\cA) = \bigoplus\limits_{[Z] \in \cI} \widetilde{\QQ}[Z]^+  \otimes_{\widetilde{\QQ}} \widetilde{\QQ}[K]
\,\, \mbox{and} \, \,
H^-(\cA) = \bigoplus\limits_{[Z] \in \cI} \widetilde{\QQ}[Z]^-  \otimes_{\widetilde{\QQ}} \widetilde{\QQ}[K].
$$
and  $H^\pm(\cA) = H(\cA)$  as $\widetilde{\QQ}$-algebras.
Let $a = [Z] K_\gamma$ and
$$
\Delta(a) = \sum\limits_{i} a_i^{(1)} \otimes a_i^{(2)} =
 \sum\limits_{[X], \, [Y] \in \cI} v^{-\langle \bar{X}, \bar{Y}\rangle}
\frac{\displaystyle P^{Z}_{X, Y}}{a_Z} [X] K_{\bar{Y} + \gamma} \otimes [Y] K_\gamma.
$$
Then we denote
$$
\Delta(a^\pm) = \sum\limits_{i} a_i^{(1)\pm} \otimes a_i^{(2)\pm} =
 \sum\limits_{[X],  \, [Y] \in \cI} v^{-\langle \bar{X}, \bar{Y}\rangle}
\frac{\displaystyle P^{Z}_{X, Y}}{a_Z} [X]^\pm K_{\pm \bar{Y} + \gamma} \otimes [Y]^\pm  K_\gamma.
$$

\begin{definition}
The Drinfeld double of $H(\cA)$ with respect to the Green's pairing $(\,-\,,\,-\,)$ is
the associative algebra $\widetilde{D}H(\cA)$, which is the free product of algebras $H^+(\cA)$ and
$H^-(\cA)$ subject to the following  relations $D(a, b)$ for all  $a, b \in H(\cA)$:
$$
\sum\limits_{i,j} a_{i}^{(1)-} b_{j}^{(2)+}  \bigl(a_{i}^{(2)}, b_{j}^{(1)}\bigr) =
\sum\limits_{i,j} b_{j}^{(1)+} a_{i}^{(2)-}  \bigl(a_{i}^{(1)}, b_{j}^{(2)}\bigr).
$$
\end{definition}

\noindent
The following proposition is well-known, see for example \cite[Section 3.2]{Joseph} for the case
of Hopf algebras and
\cite[Appendix B]{HallEll} for the case of topological bialgebras.

\begin{theorem}\label{T:DrinfDouble}
The multiplication morphism $\mathrm{mult}: H^+(\cA) \otimes_{\widetilde{\QQ}} H^-(\cA) \rightarrow \widetilde{D}H(\cA)$
is a isomorphism of $\widetilde{\QQ}$--vector spaces.  Moreover, if $\cA$ is an abelian category
 of finite length,
 then $\widetilde{D}H(\cA)$ is also a Hopf algebra such that the above morphism
 $H^+(\cA) \rightarrow \widetilde{D}H(\cA),$ $a \mapsto a \otimes \mathbbm{1}^-$ is an  injective homomorphisms of Hopf algebras.
\end{theorem}

\noindent
The following definition is due to Xiao \cite{Xiao}.
\begin{definition}
Let $\cA$ be a $\kk$--linear finitary hereditary category. The \emph{reduced} Drinfeld double
$DH(\cA)$ is the quotient of $\widetilde{D}H(\cA)$ by the two-sided ideal
$$
I = \left\langle K_{\alpha}^+ \otimes K_{-\alpha}^- - \mathbbm{1}^+ \otimes \mathbbm{1}^-
\, \, | \, \alpha \, \in K\right\rangle.$$
 Note for a finite length abelian category $\cA$,  $I$ is a Hopf ideal and
$DH(\cA)$ is also a Hopf algebra.
\end{definition}

\begin{corollary}\label{C:DrinfDoubleStr}
We have an isomorphism of $\widetilde{\QQ}$--vector spaces
$$
\mathrm{mult}: \bar{H}^{+}(\cA) \otimes_{\widetilde{\QQ}}  \widetilde{\QQ}[K] \otimes_{\widetilde{\QQ}} \bar{H}^{-}(\cA) \lar DH(\cA).
$$
\end{corollary}

\noindent
Next, we shall need
the following statement, a first version of which dates back to Dold
\cite{Dold}.

\begin{theorem}\label{T:Dold}
Let $\cA$ be a hereditary abelian category. Then any indecomposable object of the derived
category $D^b(\cA)$ is isomorphic to $X[n]$, where $X$ is an indecomposable object of $\cA$.
\end{theorem}

\begin{remark}
In what follows, we shall identify an object $X \in \mathrm{Ob}(\cA)$ with its image
under the canonical functor $\cA \rightarrow D^b(\cA)$.
\end{remark}

\noindent
The following important  theorem was recently proven by Cramer \cite{Cramer}.

\begin{theorem}\label{T:Cramer}
Let $\cA$ and $\cB$ be two $\kk$-linear finitary hereditary categories.  Assume one
 of them is artinian and there is an equivalence
of triangulated categories $D^b(\cA) \stackrel{\FF}\longrightarrow D^b(\cB)$.
Then there is an algebra
isomorphism
$$
\FF: DH(\cA) \lar DH(\cB)
$$
uniquely determined  by the following properties. For any object $X \in \mathrm{Ob}(\cA)$ such that
$\FF(X) \cong \widehat{X}[n]$ with $\widehat{X} \in \mathrm{Ob}(\cB)$ and
$n\in \ZZ$ we have:
$$
\FF\bigl([X]^\pm\bigr) = v^{-n \bigl\langle \bar{X}, \bar{X}\bigr\rangle}
\bigl[\widehat{X}\bigr]^{\pm \varepsilon(n)}
\bigl(K_{\widehat{X}}^{\pm \varepsilon(n)}\bigr)^n,
$$
where $\varepsilon(n) = (-1)^n$.
For $\alpha \in K$ we have:
$\FF(K_\alpha) = K_{\FF(\alpha)}$.
\end{theorem}

\section{Composition algebra of the Kronecker quiver}

\noindent
In this section, we study properties of  the composition algebra of the Kronecker quiver
$$
\overr{Q}
=
\xymatrix{
1  \ar@/^/[rr]^a \ar@/_/[rr]_b & & 2
}
$$
and the reduced Drinfeld double of its Hall algebra.

\begin{definition}
Consider the pair of  reflection functors $\SS^{\pm}: \Rep(\overr{Q}) \rightarrow \Rep(\overr{Q})$ defined as follows, see \cite[Section 1]{BGP}. For a representation
   $$X  = \xymatrix{
V   \ar@/^/[rr]^{A} \ar@/_/[rr]_{B} & & W}
$$
consider the short exact sequences
$$
0 \lar U' \xrightarrow{\left(\begin{smallmatrix} C' \\ D'\end{smallmatrix}\right)} V \oplus V
\xrightarrow{(A, B)} W \quad \mbox{and} \quad
V \xrightarrow{\left(\begin{smallmatrix} A \\ B \end{smallmatrix}\right)}
W \oplus W \xrightarrow{(C'', \; D'')} U'' \lar 0.
$$
Then we have:
$$
\SS^{+}(X) = \Bigl(\xymatrix{
U'   \ar@/^/[rr]^{C'} \ar@/_/[rr]_{D'} & & V}  \Bigr) \quad \mbox{and} \quad
\SS^{-}(X) = \Bigl(\xymatrix{
W   \ar@/^/[rr]^{C''} \ar@/_/[rr]_{D''} & & U''}  \Bigr).
$$
The action of $\SS^\pm$ on morphisms is defined using the universal property
of kernels and cokernels.
Note that the functor $\SS^+$ is left exact whereas $\SS^-$ is right exact.
\end{definition}

\noindent
The following theorem  summarizes main properties of the  functors
$\SS^\pm$.

\begin{theorem}\label{T:ReflFunct}
Let $\overr{Q}$ be the Kronecker quiver, $A = k\overr{Q}$ be its path algebra
 and $\cA = \Rep(\overr{Q}) = A-\mod$. Then the following properties hold:
\begin{enumerate}
\item The functors $\SS^+$ and $\SS^-$ are adjoint, i.e.~for any $X, Y \in \mathrm{Ob}(\cA)$ we
have:
$$
\Hom_\cA\bigl(\SS^-(X), Y\bigr) \cong \Hom_\cA\bigl(X, \SS^+(Y)\bigr).
$$
\item The  derived functors $\RR\SS^+$ and $\LL\SS^-$ are also adjoint. Moreover, they are
mutually inverse auto-equivalences
of the derived category $D^b(\cA)$.
\item Let $X \in \mathrm{Ob}(\cA)$ be an indecomposable object. Then we have:
$$
\RR^{0} \SS^{+}(X) =
\left\{
\begin{array}{ccc}
\SS^+(X)   & \, \mbox{\rm if}\, & \, X \not\cong S_2 \\
0  & \, \mbox{\rm if}\, & \, X \cong S_2
\end{array}
\right.
\quad \mbox{\rm and} \quad \RR^{1} \SS^{+}(X) =
\left\{
\begin{array}{ccc}
0  & \, \mbox{\rm if}\, & \, X \not\cong S_2\\
S_1 & \, \mbox{\rm if}\, & \, X \cong S_2.
\end{array}
\right.
$$
Similarly, we have:
$$
\LL^{0} \SS^{-}(Y) =
\left\{
\begin{array}{ccc}
\SS^-(Y)   & \, \mbox{\rm if}\, & \, Y \not\cong S_1 \\
0  & \, \mbox{\rm if}\, & \, Y \cong S_1
\end{array}
\right.
\quad \mbox{\rm and} \quad
\LL^{-1} \SS^{-}(Y) =
\left\{
\begin{array}{ccc}
0  & \, \mbox{\rm if}\, & \, Y \not\cong S_1 \\
S_2 & \, \mbox{\rm if}\, & \, Y \cong S_1.
\end{array}
\right.
$$
\item
In particular, the reflection functors $\SS^{-}$ and $\SS^{+}$ yield
  mutually inverse equivalences between the categories $\Rep\bigl(\overr{Q}\bigr)^1$ and $\Rep\bigl(\overr{Q}\bigr)^2$, which
are the full subcategories of $\Rep\bigl(\overr{Q}\bigr)$
consisting of objects without direct summands isomorphic to $S_1$ and $S_2$ respectively.
\item Let $\nu = \DD\bigl(\Hom_A(\,-\,,A)\bigr): A-\mod \rightarrow A-\mod$ be the Nakayama
functor and $\SS:= \LL\nu: D^b(\cA) \rightarrow D^b(\cA)$ be its derived functor. Then we have an  isomorphism of vector spaces
$$
\Hom_{D^b(\cA)}\bigl(X, \,  \SS(Y)\bigr) \lar  \DD \Hom_{D^b(\cA)}(Y,\,  X),
$$
functorial in both arguments, where $\DD$ is the duality over $\kk$.
 In other words, $\SS$ is the Serre functor
of  the triangulated category $D^b(\cA)$ in the sense of \cite{BondalKapranov}.
\item The  functors $\SS$  and   $\SS^+$ are related by an isomorphism:
$
\SS \cong (\RR\SS^+)^2[1].
$
\end{enumerate}
\end{theorem}

\begin{proof}
The first part was essentially proven in \cite[Section 1]{BGP}. There the  authors construct
two natural transformations of functors $\imath: \SS^- \SS^+ \rightarrow \mathbbm{1}$
and $\jmath: \mathbbm{1} \rightarrow \SS^+ \SS^-$. It can be easily shown that
they define mutually inverse bijections $$\Hom_\cA\bigl(\SS^-(X), Y\bigr)
\longleftrightarrow
\Hom_\cA\bigl(X, \SS^+(Y)\bigr).$$ See also \cite[Section VII.5]{ASS} for a proof using
tilting functors.

The fact that the derived functors  $\RR \SS^+$ and $\LL \SS^-$ are adjoint, is a general property of
an adjoint pair, see for example \cite[Lemma 15.6]{Keller}. For the proof that
$\RR \SS^+$ and $\LL \SS^-$ are equivalences of categories, see for example \cite[Section VII.5]{ASS}.

By Theorem \ref{T:Dold}, the complexes  $\RR \SS^+(X)$ and $\LL \SS^-(X)$ have exactly one non-vanishing
cohomology for  an indecomposable object $X \in \mathrm{Ob}(\cA)$. This proves the formulae
listed in the third item. The fourth statement is proven in \cite[Section 1]{BGP} and for
 the fifth we refer to  \cite[Section 4.6]{Happel}.

For a proof of the last statement,  first note that the Auslander-Reiten functor
$$\tau = \DD \Ext^1_A(\,-\,,A): \; \Rep(\overr{Q}) \longrightarrow \Rep(\overr{Q})$$ is isomorphic to
the Coxeter functor  $\AA^+ := (\SS^+)^2$, see  \cite[Section 5.3]{Gabriel} and
 \cite[Proposition II.3.2]{Tepetla}. Next,  the canonical transformation of functors
 $\RR \AA^+ \rightarrow (\RR \SS^+)^2$ is an isomorphism on the indecomposable
 injective modules $I(1)$ and $I(2)$, hence it is an isomorphism. Finally, by
 \cite[Proposition I.7.4]{RD} we know that the derived functors
 $\RR \tau[1]$ and $\LL\nu$ are isomorphic.

\end{proof}

\begin{remark}\label{R:SerreFexplicit}
Let $i \in Q_0 = \left\{1,2\right\}$ be a vertex and $P(1) = A e_i$ be the  indecomposable
projective module, which is the projective cover of the simple module $S_i$. Then
$\LL\nu\bigl(P(i)\bigr) = \nu\bigl(P(i)\bigr) = I(i)$, where $I(i)$ is the injective envelope of $S_i$.
\end{remark}

\noindent
Let $X$ be an indecomposable object of $\cA$ and  $B = \End_{\cA}(X)$.
By  Serre duality, we have a canonical isomorphism of  $B$-bimodules
$$\Hom_\cA(X, X)  \lar   \DD\Bigl(\Hom_{D^b(\cA)}\bigl(X, \SS(X)\bigr)\Bigr).$$
 Let $w$ be a non-zero element
 of the  socle of $\Hom_{D^b(\cA)}\bigl(X, \SS(X)\bigr)$ viewed as the right  $B$-module.

\begin{lemma}\label{L:LemmaHappel}
Consider a  distinguished triangle $$\SS[-1](X) \stackrel{u}\lar Y \stackrel{v}\lar X \stackrel{w}\lar \SS(X)$$  given by the morphism $w$. Then this triangle is almost split.
Moreover, if $X$ is non-projective,
then $H^i\bigl(\SS(X)\bigr) = 0$ for $i \ne 1$ and  the short exact sequence
$$
0 \lar \tau(X) \stackrel{u}\lar Y \stackrel{v}\lar X  \lar 0
$$
is almost split, where $\tau(X) = \DD \Ext^1_A(X, A) \cong H^1\bigl(\SS(X)\bigr)$.
\end{lemma}

\begin{proof}
For the first part of the statement, see the proof of \cite[Proposition I.2.3]{ReitenBergh}.
For the second, see \cite[Section 4.7]{Happel}.
\end{proof}

\begin{definition}
 An
object $X \in \mathrm{Ob}(\cA)$ is called
\begin{enumerate}
\item  pre-projective if there exists a projective object $P$ and $m \ge 0 $ such that $X \cong \tau^{-m}(P)$,
\item pre-injective if there exists an injective object $I$ and $m \ge 0 $ such that $X \cong \tau^{m}(I)$.
\end{enumerate}
\end{definition}

\noindent
Recall the classification of the indecomposable objects of the category of representations
of $\overr{Q}$ over an arbitrary  field $\kk$, see for example \cite[Section 3.2]{RingelBook}.

\begin{theorem}\label{T:indofKron}
The indecomposable representations of the Kronecker quiver $\overr{Q}$ over an arbitrary field
$\kk$ are the following.
\begin{enumerate}
\item Indecomposable pre-projective objects $$P_n = \xymatrix{
\kk^n   \ar@/^/[rr]^{A^{\mathrm{pro}}_n} \ar@/_/[rr]_{B^{\mathrm{pro}}_n} & & \kk^{n+1}
}, \quad n \ge 0,$$ where
$$
A^{\mathrm{pro}}_n =
\left(
\begin{array}{cccc}
1      & 0      &  \dots    &   0 \\
0      & 1      &  \dots      &   0 \\
\vdots & \vdots & \ddots &   \vdots          \\
0      & 0      & \dots  & 1       \\
0      & 0      & \dots  & 0
\end{array}
\right) \quad \mbox{and}
\quad
B^{\mathrm{pro}}_n =
\left(
\begin{array}{cccc}
0      & 0      & \dots        & 0 \\
1    &  0     &  \dots     &   0 \\
0      & 1      &  \dots     &   0 \\
\vdots & \vdots & \ddots &   \vdots          \\
0      & 0      & \dots  &  1
\end{array}
\right).
$$
In particular, $P_0$ and $P_1$ are the indecomposable projective objects.
\item Indecomposable pre-injective objects
$$I_n = \xymatrix{
\kk^{n+1}   \ar@/^/[rr]^{A^{\mathrm{inj}}_n} \ar@/_/[rr]_{B^{\mathrm{inj}}_n} & & \kk^{n}
}, \quad n \ge 0,$$ where
$$
A^{\mathrm{inj}}_n =
\left(
\begin{array}{ccccc}
1      & 0      &  0     &  \dots & 0 \\
0      & 1      &  0     &  \dots & 0 \\
\vdots & \vdots & \ddots &  \ddots & \vdots          \\
0      & 0      & \dots  & 1      & 0
\end{array}
\right) \quad \mbox{and}
\quad
B^{\mathrm{inj}}_n =
\left(
\begin{array}{ccccc}
0     & 1      &  0     &  \dots & 0 \\
0      & 0      &  1     &  \dots & 0 \\
\vdots & \vdots & \ddots &  \ddots & \vdots          \\
0      & 0      & \dots  & 0      & 1
\end{array}
\right).
$$
In particular, $I_0$ and $I_1$ are the indecomposable injective  objects.
\item Tubes $$T_{n, \pi} = \xymatrix{
\kk^{nl}   \ar@/^/[rr]^{A^{\mathrm{tub}}_{n, \pi}} \ar@/_/[rr]_{B^{\mathrm{tub}}_{n, \pi}} & & \kk^{nl}
}, \quad n \ge 1,$$
where $\pi(x, y) \in \kk[x,y]$ is an irreducible homogeneous polynomial of degree $l$.
For $(\pi) \ne (x)$  we have:  $A^{\mathrm{tub}}_{n, \pi} = I_{nl}$ and $B^{\mathrm{tub}}_{n, \pi} = F\bigl(\pi(1,y)^n\bigr)$ is the   Frobenius normal form defined by  the polynomial  $\pi(1, y)^n$.
For $(\pi) =  (x)$  we set $A^{\mathrm{tub}}_{n, \pi} = F(x^n)$ and $B^{\mathrm{tub}}_{n, \pi} = I_n$.
\end{enumerate}
 Moreover, the  Auslander--Reiten quiver of $\Rep(\overr{Q})$ has the following form:
$$
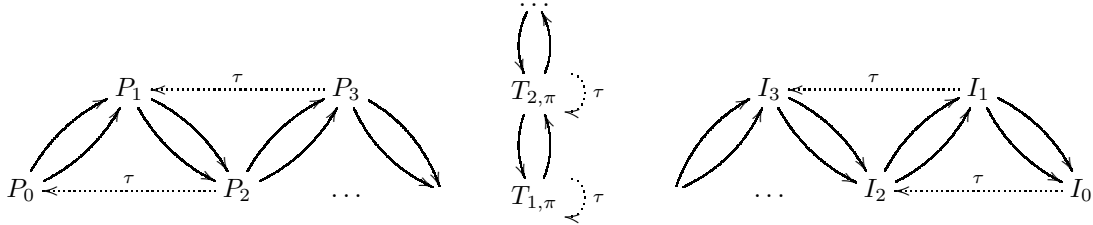
\xymatrix
{ & & &  & \\
                               & P_1 \ar@/^/[dr]  \ar@/_/[dr] &                               & \ar@{.>}[ll]_\tau P_3   \ar@/^/[dr]  \ar@/_/[dr] & \\
 P_0  \ar@/^/[ur]  \ar@/_/[ur] &                               & P_2 \ar@/^/[ur]  \ar@/_/[ur]  \ar@{.>}[ll]_\tau & \dots &
}
\qquad
\xymatrix
{  \dots \ar@/_/[d] \\
T_{2, \pi} \ar@/_/[d] \ar@/_/[u] \ar@{.>}@(ur,rd)^{\tau} \\
T_{1, \pi} \ar@/_/[u] \ar@{.>}@(ur,rd)^{\tau}
 }
 \qquad
\xymatrix
{
& & &  & \\
                               & I_3 \ar@/^/[dr]  \ar@/_/[dr] &                               & \ar@{.>}[ll]_\tau I_1   \ar@/^/[dr]  \ar@/_/[dr] & \\
    \ar@/^/[ur]  \ar@/_/[ur] &       \dots                         & I_2 \ar@/^/[ur]  \ar@/_/[ur]    &  &  \ar@{.>}[ll]_\tau I_0
}
$$
 whereas the Auslander--Reiten quiver
of the derived category  $D^b\bigl(\Rep(\overr{Q})\bigr)$ is the  union of components obtained by applying the shift functor to
$$
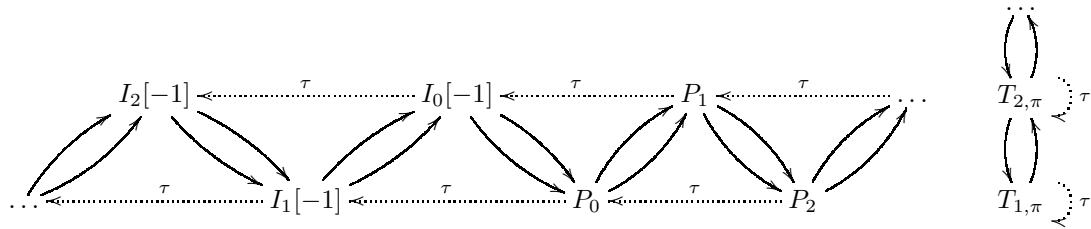
\xymatrix
{
&  & &  &
 & &   & \\
 & I_2[-1] \ar@/^/[dr]  \ar@/_/[dr] & & \ar@{.>}[ll]_\tau I_0[-1]   \ar@/^/[dr]  \ar@/_/[dr] &
 & \ar@{.>}[ll]_\tau P_1   \ar@/^/[dr]  \ar@/_/[dr] &  & \ar@{.>}[ll]_\tau \dots \\
 \dots \ar@/^/[ur]  \ar@/_/[ur] & &  I_1[-1] \ar@/^/[ur]  \ar@/_/[ur]  \ar@{.>}[ll]_\tau
 & & P_0  \ar@/^/[ur]  \ar@/_/[ur]  \ar@{.>}[ll]_\tau & &
 P_2 \ar@/^/[ur]  \ar@/_/[ur]  \ar@{.>}[ll]_\tau & \\
}
\qquad
\xymatrix
{  \dots \ar@/_/[d] \\
T_{2, \pi} \ar@/_/[d] \ar@/_/[u] \ar@{.>}@(ur,rd)^{\tau}  \\
T_{1, \pi} \ar@/_/[u] \ar@{.>}@(ur,rd)^{\tau}
 }
$$
\end{theorem}

\medskip
\noindent
Now we return to our study of Hall algebras of quiver.
Applying Theorem \ref{T:Cramer}, we get the following corollary, which is due to Sevenhant
and van den Bergh \cite{SevenhBergh}, see also \cite{XiaoYang}.

\begin{corollary}
The derived reflection functor $\RR \SS^{+}$  induces an algebra isomorphism of the
 Drinfeld doubles:
$\SS^{+}: DH(\overr{Q}) \rightarrow  DH(\overr{Q}),
$
whose inverse is  induced by the adjoint functor
$\LL \SS^{-}$.
\end{corollary}

\begin{definition}
Let $S_1$ and $S_2$ be the  simple objects of $\cA$. Consider the subalgebra $\bar{C}(\overr{Q})$  of the Hall algebra $H(\cA)$ generated
by $[S_1]$ and  $[S_2]$. The composition subalgebra is
defined as follows:
$$C(\overr{Q}) :=
\bar{C}(\overr{Q}) \otimes_{\widetilde{\QQ}} \widetilde{\QQ}[K].$$ Note that
$C(\overr{Q})$ is a Hopf subalgebra of $H(\cA)$.
\end{definition}

\noindent
Our next goal is to show the automorphisms  $\SS^{\pm}$  of $DH(\overr{Q})$ map the
reduced Drinfeld double
of the composition algebra $C(\overr{Q})$ to itself. In other words, one has to check that both  simple modules  $S_i \in \mathrm{Ob}\bigl(\cA\bigr)$, $i = 1,2$  we have:
$\SS^{\pm}\bigl([S_i]\bigr) \in DC(\overr{Q})$.
Note that in the \emph{notations
of Theorem \ref{T:indofKron}} we have: $\SS^+(S_1) = I_1$
is an  indecomposable injective  module and $\SS^-(S_2) = P_1$ is an  indecomposable projective
 module. Hence, it is sufficient to check the following lemma.

\begin{lemma}\label{P:refloncomp}
The elements $[P_1]$ and $[I_1]$ belong to the composition algebra $C(\overr{Q})$.
\end{lemma}

\begin{proof}
 Using a straightforward calculation, we get
 the following explicit formulae  for the
classes of the non-simple  indecomposable projective and injective  modules:
$$
[P_1] = \sum\limits_{a+b=2} (-1)^a v^{-b} [S_2]^{(a)} \circ [S_1] \circ [S_2]^{(b)}
\qquad
[I_1] = \sum\limits_{a+b=2} (-1)^a v^{-b} [S_1]^{(a)} \circ [S_2] \circ [S_1]^{(b)},
$$
where $[X]^{(n)} = \frac{\displaystyle [X]^n}{\displaystyle [n]_v!}$ for an object $X$ of the category $\cA$.
\end{proof}

\begin{definition}\label{D:quantaff}
Let $C = \left(\begin{smallmatrix}2 & -2 \\ -2 & 2\end{smallmatrix}\right)$
be the Cartan matrix of the affine Lie algebra
$\mathfrak{g} = \widehat{\mathfrak{sl}}_2$.
The Hopf algebra $U_v\bigl(\widehat{\mathfrak{sl}}_2\bigr)$ is generated over
$R$ by the elements
$E_1, E_2$, $F_1, F_2$, $K_1^\pm, K_2^\pm$ subject to the relations

\begin{enumerate}
\item the elements $Z^\pm  := K_1^\pm  K_2^\pm$ are  central;
\item $K_i^\pm K_i^\mp = 1 = K_i^\mp K_i^\pm$, $\, i = 1,2$;
\item  $K_i  E_j = v^{-c_{ij}} E_j K_i$ and   $K_i  F_j = v^{c_{ij}} F_j K_i$, $\, i = 1,2$;
\item $[E_i, F_j] = \delta_{ij} v \frac{\displaystyle K_i - K_i^{-1}}{\displaystyle v - v^{-1}}$,
     $\, i = 1,2$;
\item $\sum_{k=0}^{3} (-1)^k E_i^{(k)} E_j  E_i^{(3-k)} = 0$  for $1 \le i \ne j \le 2$;
\item $\sum_{k=0}^{3} (-1)^k F_i^{(k)} F_j  F_i^{(3-k)} = 0$ for $1 \le i \ne j \le 2$.
\end{enumerate}

\noindent
The Hopf algebra structure is given by the following formulae:
\begin{enumerate}
\item $\Delta(E_i) = E_i \otimes \mathbbm{1} + K_i \otimes E_i$,  $\Delta(F_i) = F_i \otimes K_i^{-1} + \mathbbm{1} \otimes F_i$ and  $\Delta(K_i) = K_i \otimes K_i$;
\item $\eta(E_i) = 0 = \eta(F_i)$, $\eta(K_i) = 1$;
\item $S(E_i) = - K_i^{-1} E_i$, $S(F_i) = -F_i K_i$, $S(K_i) = K_i^{-1}$ for all $i = 1,2$.
\end{enumerate}
\end{definition}

\medskip
\noindent
The following result is a special  case of a more general statement, which is essentially  due to Ringel
\cite{Ringel} and Green \cite{Green}.

\begin{theorem}\label{T:RingelGreen}
The $\widetilde\QQ$--linear
morphism $U_q(\mathfrak{g}) :=  U_v(\mathfrak{g}) \otimes_R \widetilde{\QQ}
\stackrel{\mathsf{ev}_q}\lar  DC(\overr{Q})$ mapping $E_i$ to $[S_i]^+$, $F_i$ to $[S_i]^-$
and $K_i$ to $K_{\bar{S}_i}$ for $i = 1,2$, is an isomorphism of Hopf algebras.
Moreover, if we set
$$DC_{\mathrm{gen}}(\overr{Q}):= \prod\limits_{q \in \kP} DC\bigl(\Rep(\FF_q \overr{Q})\bigr)$$
then the $R$--linear  map $\mathsf{ev} = \prod{\mathsf{ev}_q}: U_v (\mathfrak{g}) \rightarrow DC_{\mathrm{gen}}(\overr{Q})$
is injective. The same applies  to the subalgebra $U_v(\mathfrak{g}^+ \oplus \mathfrak{h})$ and
the algebra $C_{\mathrm{gen}}(\overr{Q}):= \prod_{q \in \kP} C\bigl(\Rep(\FF_q \overr{Q})\bigr)$.
\end{theorem}

\begin{corollary}\label{C:ReflonCompAlg}
The derived functors $\RR \SS^+$ and $\LL \SS^-$  induce a pair  of mutually inverse
 automorphisms  $\SS^{\pm}$
of the algebra
 $U_v(\mathfrak{g})$ such that the following diagrams are  commutative:
$$
\xymatrix
{
U_v(\mathfrak{g}) \ar[rr]^{\SS^{+}} \ar[d]_{\ev_q} & & U_v(\mathfrak{g}) \ar[d]^{\ev_q} \\
DC(\overr{Q}) \ar[rr]^{\SS^{+}} & & DC(\overr{Q})
}
\qquad
\xymatrix
{
U_v(\mathfrak{g}) \ar[rr]^{\SS^{-}} \ar[d]_{\ev_q} & & U_v(\mathfrak{g}) \ar[d]^{\ev_q} \\
DC(\overr{Q}) \ar[rr]^{\SS^{-}} & & DC(\overr{Q}).
}
$$
Their  action on the generators is   given by the  following formulae:
$$
\begin{array}{|l||l|}
\hline
E_1  \xrightarrow{\SS^{+}}       \sum\limits_{a+b \,=\, 2}
(-1)^a v^{-b} E_1^{(a)}  E_2  E_1^{(b)}     &
E_2  \xrightarrow{\SS^{-}}       \sum\limits_{a+b \,=\, 2}
(-1)^a v^{-b} E_2^{(a)}  E_1  E_2^{(b)}
 \\
\hline
F_1  \xrightarrow{\SS^{+}}       \sum\limits_{a+b \,=\, 2}
(-1)^a v^{-b} F_1^{(a)}  F_2  F_1^{(b)}     &
F_2  \xrightarrow{\SS^{-}}       \sum\limits_{a+b \,=\, 2}
(-1)^a v^{-b} F_2^{(a)}  F_1  F_2^{(b)}
 \\
 \hline
 E_2  \xrightarrow{\SS^{+}}  v^{-1}  F_1 K_1, \quad    F_2 \xrightarrow{\SS^{+}}  v^{-1}
E_1 K_1^{-1}  &
E_1  \xrightarrow{\SS^{-}}  v F_2 K_2^{-1}, \quad    F_1 \xrightarrow{\SS^{-}}  v  E_2 K_2 \\
 \hline
K_1  \xrightarrow{\SS^{+}} K_1^2 K_2, \quad  K_2  \xrightarrow{\SS^{+}} K_1^{-1} &
 K_1 \xrightarrow{\SS^{-}} K_2^{-1},  \quad K_2  \xrightarrow{\SS^{-}} K_1 K_2^2   \\
\hline
\end{array}
$$
As it was explained in \cite[Theorem 13.1]{SevenhBergh},
in  the conventions of  Remark \ref{R:formofrelat}, these automorphisms  are
the symmetries discovered by Lusztig \cite{Lusztig} (strictly speaking,
Lusztig's symmetries are obtained by composing  the automorphisms $\SS^\pm$ with
a flip  sending $E_1$ to $E_2$, $F_1$ to $F_2$ and $K_1$ to $K_2$. This happens
because we wish to consider the functors $\SS^\pm$ as \emph{endofunctors} of $\Rep(\overr{Q})$
whereas in \cite{SevenhBergh} one interchanges the direction of the arrows).
\end{corollary}

\begin{remark}\label{R:formofrelat}
Note that the relations of the quantum affine algebra $U_v(\widehat{\mathfrak{sl}}_2)$
given in Definition \ref{D:quantaff},   slightly differ from the classical ones as defined for instance in \cite{ChariPressley, Lusztig, JimboMiwa}.
Namely, we impose the commutation relation
$[E_i, F_i] = v \frac{\displaystyle K_i - K_i^{-1}}{\displaystyle v - v^{-1}}$,
whereas the  conventional form   would be
 $[{E}_i, {F}_i] =  \frac{\displaystyle K_i - K_i^{-1}}{\displaystyle v - v^{-1}}$, $\, i = 1,2$.
 However, one can easily pass to the conventional form replacing at the \emph{first step}
 $v$ by $v^{-1}$ and then applying  the Hopf algebra automorphism
 sending $E_i$ to $E_i$, $K_i$ to $K_i$ and $F_i$ to $-v^{-1} F_i$ for $i = 1,2$ at the \emph{second step}.

 \vspace{2mm}
 \noindent
 In order to get the relations of the Drinfeld double, which are closer  to the
  conventional ones, one can alternatively redefine Green's form by setting
 $$
\Bigl([X] K_\alpha, [Y] K_\beta\Bigr)_{\mathrm{new}} = v^{-(\alpha, \beta)- \dim(X)}
\frac{\displaystyle \delta_{X, Y}}{\displaystyle a_X},
$$
where $\dim(X)$ is the dimension of $X$ over $\kk$ viewed as an $\kk\overr{Q}$--module, and
take the reduced Drinfeld double with respect to $(\,-\,,\,-\,)_{\mathrm{new}}$.

\vspace{2mm}
\noindent
Nevertheless, we  prefer
 to follow the form of the relations as stated  in Definition \ref{D:quantaff}, because they seem
 to be more natural from the point of view of Hall algebras.
\end{remark}

\noindent
The following automorphism  plays an important role in our study of the
algebra $U_v\bigl(\widehat{\mathfrak{sl}}_2\bigr)$.

\begin{definition}\label{D:CoxAutonComp}
The automorphism
$\AA = \bigl(\SS^+\bigr)^2  :  DC(\overr{Q}) \rightarrow  DC(\overr{Q})$
is called \emph{Coxeter automorphism} of $DC(\overr{Q})$. Using Corollary \ref{C:ReflonCompAlg}, we
also  obtain  the corresponding  automorphism of the algebra $U_v(\mathfrak{g})$, given by the commutative diagram
$$
\xymatrix
{
U_v(\mathfrak{g}) \ar[rr]^{\AA^+} \ar[d]_\ev & & U_v(\mathfrak{g})  \ar[d]^\ev \\
DC_{\mathrm{gen}}(\overr{Q}) \ar[rr]^{(\SS^+)^2} &   & DC_{\mathrm{gen}}(\overr{Q}).
}
$$
The inverse automorphism $\AA^- = \bigl(\SS^-\bigr)^2$ is defined in a similar way.
\end{definition}

\noindent
Using this technique, we can give a new proof of  the following well-known result, see \cite{Zhang, Szanto}.

\begin{lemma}
The composition algebra $C(\overr{Q})$ contains all the elements  $[X]$, where
$X$ is either an indecomposable pre-injective module or an indecomposable pre-pro\-jecti\-ve
module.
\end{lemma}

\begin{proof} Let $X$ be either pre-projective or pre-injective indecomposable representation of
$\overr{Q}$.
Since the automorphism $\AA^+$ acts on the Drinfeld double $DC(\overr{Q})$, we know that
$\AA^+\bigl([X]\bigr)$ belongs to $DC(\overr{Q})$.  Moreover, $\AA^+\bigl([X]\bigr)$
has the form $[Y]^\pm K_\alpha$ for some
indecomposable pre-projective or pre-injective representation $Y$ and some $\alpha \in K$.
As in Corollary
\ref{C:DrinfDoubleStr}, we have a triangular decomposition
 $$DC(\overr{Q}) = \bar{C}(\overr{Q})^+ \otimes_{\widetilde\QQ}
\widetilde{\QQ}[K] \otimes_{\widetilde\QQ} \bar{C}(\overr{Q})^-.$$
From this fact it follows that
$
DC(\overr{Q}) \cap {H}(\overr{Q})^\pm   = {C}(\overr{Q})^\pm := \bar{C}(\overr{Q})^+ \otimes_{\widetilde\QQ}
\widetilde{\QQ}[K].
$
 Hence, the element
$\AA^+\bigl([X]\bigr)$ has to  belong to one of the aisles
${C}(\overr{Q})^\pm$.
The proof for $\AA^-\bigl([X]\bigr)$ is analogous.
Since any indecomposable pre-projective  or pre-injective object
$X$ is (up to a shift) of the form $\AA^{m}(S)$, where $m \in \ZZ$ and $S$ is a simple object in
$\Rep(\overr{Q})$, this implies the claim.
\end{proof}

\section{Composition algebra of the category of coherent sheaves on  $\PP^1$}

\noindent
In this subsection, we consider the composition subalgebra of the category of coherent
sheaves on the  projective line $\PP^1$.

First note that  the maps $\Pic(\PP^1)  \xrightarrow{\deg} \ZZ$ and $K_0\bigl(\Coh(\PP^1)\bigr)
\xrightarrow{(\rk, \deg)} \ZZ^2$ are isomorphisms of abelian groups. Next,  recall
some well-known facts on coherent sheaves on $\PP^1$.

\begin{theorem}\label{T:basicsonP1}
The indecomposable objects of the category  $\Coh(\PP^1)$ are:
\begin{enumerate}
\item  line bundles $\kO_{\PP^1}(n)$, $n \in \ZZ$;
\item torsion sheaves  $\kT_{t,x} := \kO_{\PP^1}/\mathfrak{m}_x^t$, where $x \in \PP^1$ is a closed
point and $t \in \ZZ_{>0}$.
\end{enumerate}
\end{theorem}

\begin{definition}\label{D:compalggener}
Let $\Tor(\PP^1)$ be the abelian category of torsion coherent sheaves on $\PP^1$ and $H(\PP^1)_{\mathrm{tor}} \subseteq H(\PP^1)$ be its Hall algebra.
For any integer $r \ge 1$ consider the element
$$
\mathbbm{1}_{(0, r)} := \sum\limits_{\kT \in \Tor(\PP^1): \, \bar{\kT} = (0, r)} [\kT] \in
 H(\PP^1)_{\mathrm{tor}}.
$$
Next, consider the family  $\bigl\{T_r\bigr\}_{r \ge 1}$  determined
 by  $\bigl\{\mathbbm{1}_{(0,r)}\bigr\}_{r \ge 1}$ using
 the generating series
$$
1 + \sum\limits_{r = 1}^\infty \mathbbm{1}_{(0, r)} t^r = \exp\Bigl(\sum\limits_{r = 1}^\infty
\frac{T_r}{[r]_v} t^r\Bigr).
$$
Finally, the elements $\bigl\{\Theta_r\bigr\}_{r \ge 1}$ are defined by the generating series
$$
1 + \sum\limits_{r = 1}^\infty \Theta_r  t^r = \exp\bigl((v^{-1}-v)\sum\limits_{r = 1}^\infty T_r t^r\bigr).
$$
In what follows, we  set $\mathbbm{1}_{(0, 0)} = T_0 = \Theta_0 = [0] = \mathbbm{1}$.
\end{definition}

\begin{proposition}
In the notations as above we have:
\begin{itemize}
\item Any of three families  $\bigl\{\mathbbm{1}_{(0, r)}\bigr\}_{r \ge 1}$,
$\bigl\{ T_r \bigr\}_{r \ge 1}$ and $\bigl\{ \Theta_r \bigr\}_{r \ge 1}$
introduced in Definition \ref{D:compalggener}, generates the
same subalgebra $U(\PP^1)_{\mathrm{tor}}$ of  the Hall algebra
$H(\PP^1)_{\mathrm{tor}}$;
\item For any $r,s  \ge 1$ we have the equalities:
$$
\Delta(T_r) = T_r \otimes \mathbbm{1} + K_{(0,r)} \otimes T_r \quad
\mbox{and} \quad \bigl(T_r, \, T_s\bigr) = \delta_{r,s}
\frac{[2r]}{r(v^{-1}-v)}.
$$
\end{itemize}
\end{proposition}

\begin{proof}
The first part of this proposition is trivial, a proof of the second can be found
 in \cite{OlivierNotes}.
\end{proof}

\noindent
Using this proposition, we get the following statement.

\begin{lemma}\label{L:someGreenprod} In the notation as above we have: $\left(\Theta_r, T_r\right) =
\frac{\displaystyle [2r]}{\displaystyle r}$ for any $r \in \ZZ_{>0}$.
\end{lemma}

\begin{proof}
For a sequence of non-negative integers $\underline{c} = (c_r)_{r \in \ZZ_{> 0}}$ such that
all but finitely many entries are zero, we set:
$$T_{\underline{c}} := \prod\limits_{r = 1}^\infty \frac{T_r^{c_r}}{c_r!}
\quad \mbox{and} \quad  c = \sum\limits_{r = 1}^\infty  r c_r.
$$
Then we have:
$$
\Delta(T_{\underline{c}}) = \sum\limits_{\underline{a} + \underline{b} = \underline{c}}
T_{\underline{a}} K_{(0, b)} \otimes T_{\underline{b}}.
$$
In particular, by induction we obtain:
$$
\left(T_c, T_{\underline{c}}\right) =
\left\{
\begin{array}{cl}
\frac{\displaystyle [2c]}{\displaystyle c(v^{-1}-v)} & \mbox{if} \quad \underline{c} = (0,\dots, 0,\underset{c \mathrm{-th \, pl}}{1},0,\dots), \\
0 & \mbox{otherwise}.
\end{array}
\right.
$$
Using the formula $\Theta_r = (v^{-1}- v) T_r  \, +  \mbox{ \,monomials of length\, }
\ge 2\mbox{\,\, in\,\,} T_s,$
 the statement follows.
\end{proof}

\noindent
Summing up, the first family  of generators $\bigl\{\mathbbm{1}_{(0, r)}\bigr\}_{r \ge 0}$ of the algebra $U(\PP^1)_{\mathrm{tor}}$
 has a clear algebro--geometric meaning. The second set $\bigl\{T_r\bigr\}_{r \ge 0}$
 has a good behavior with respect to the bialgebra structure:
all the generators  $T_r$ are primitive and orthogonal with respect to  Green's pairing.
The role of the third family $\bigl\{\Theta_r\bigr\}_{r \ge 0}$ is explained by the following proposition.

\begin{proposition}\label{P:comultoflb}
For any $n \in \ZZ$ we have:
$$
\Delta\bigl(\bigl[\kO_{\PP^1}(n)\bigr]\bigr) =
\bigl[\kO_{\PP^1}(n)\bigr] \otimes \mathbbm{1} + \sum\limits_{r = 0}^\infty \Theta_r K_{(1, n-r)}
\otimes \bigl[\kO_{\PP^1}(n-r)\bigr].
$$
\end{proposition}

\begin{proof}
We refer to \cite[Theorem 3.3]{Ka} or
\cite[Section 12.2]{SchiffmannInvent}  for a proof of this
result.
\end{proof}

\begin{remark}\label{R:anotherformofTheta}
In the next section, we shall need another description  of the elements $\Theta_r$, see for example
\cite[Example 4.12]{OlivierNotes}
$$
\Theta_r =  v^{-r}\sum\limits_{\substack{x_1, \dots, x_m \in \PP^1; \,x_i \ne x_j \,1 \le i \ne j \le m
\\
t_1, \dots, t_m \,:\,
\sum_{i=1}^m t_i \deg(x_i) = r}}
\prod_{i = 1}^m \bigl(1 - v^{2\deg(x_i)}\bigr)\bigl[\kT_{t_i, x_i}\bigr].
$$
\end{remark}

\begin{definition}\label{D:compalg}
The  composition algebra $U(\PP^1)$ is the subalgebra of the Hall algebra
$H\bigl(\Coh(\PP^1)\bigr)$ generated by the elements $L_n:= \bigl[\kO_{\PP^1}(n)\bigr]$, $T_r$ and $K_\alpha$,
where $n \in \ZZ$,  $r \ge 1$ and $\alpha \in K_0\bigl(\Coh(\PP^1)\bigr) \cong  \ZZ^2$.
We also use the notations: $\delta = (0,1) \in K_0(\PP^1)$, $C = K_\delta$ and $K = K_{(1,0)}$.
\end{definition}

\vspace{2mm}

\noindent
A complete list of relations between the generators of
 the composition algebra $U(\PP^1)$ was obtained  by Kapranov
\cite{Ka} and Baumann--Kassel \cite{BK}, see also \cite[Section 4.3]{OlivierNotes}.

\begin{theorem}\label{T:relationsinP1}
The elements  $L_n, T_r$, $K$ and $C$
satisfy  the following relations:
\begin{enumerate}
\item $C$ is central;
\item $[K, T_n] = 0 = [T_n, T_m]$ for all $m, n  \in \ZZ_{> 0}$;
\item $K L_n = v^{-2} L_n K$ for all $n \in \ZZ$;
\item $\bigl[T_r, L_n\bigr] = \frac{\displaystyle [2r]}{\displaystyle r} L_{n+r}$ for all $n \in \ZZ$ and $r \in \ZZ_{> 0}$;
\item $L_{m} L_{n+1}  + L_{n}L_{m+1} =  v^2 \bigl(L_{n+1} L_{m} + L_{m+1} L_{n})$ for all $m, n \in
\ZZ$.
\end{enumerate}

\vspace{1mm}
\noindent
Let $U(\PP^1)_{\mathrm{vec}}$ be the subalgebra of $U(\PP^1)$ generated by the elements
$L_n$, $n \in \mathbb{Z}$.
Then the $\widetilde\QQ$--linear map $U(\PP^1)_{\mathrm{vec}} \otimes_{\widetilde{\QQ}} U(\PP^1)_{\mathrm{tor}}
\otimes_{\widetilde{\QQ}} \widetilde{\QQ}[K]
\xrightarrow{\mathsf{mult}} U(\PP^1)$
is an isomorphism. In particular, the elements $$B_{\underline{m}, \, \underline{l}, \,  a, \, b} =
\prod_{n \in \ZZ} L_n^{m_n} \circ \prod_{r \in \ZZ^+}  T_r^{l_r} \circ K^a C^b,$$
where $a, b \in \ZZ$,
$\underline{m} = (m_n)_{n \in \ZZ}$ and $\underline{l} = (l_r)_{r \in \ZZ_{>0}}$ are
sequences of non-negative integers such that all but finitely entries are zero, form a basis
of $U(\PP^1)$.
\end{theorem}

\vspace{1mm}
\noindent
It turns out that in order to relate the reduced Drinfeld double
 $DU(\PP^1)$ with Drinfeld's  new presentation of
$U_v\bigl(\widehat{\mathfrak{sl}}_2\bigr)$ \cite{Drinfeld},
one has to  modify the definition of the generators
$T_r^\pm$ and $\Theta_r^\pm$. Similarly to Corollary \ref{C:DrinfDoubleStr},  we have a triangular decomposition
$$
DU(\PP^1) = \bar{U}(\PP^1)^+ \otimes_{\widetilde\QQ}
\widetilde\QQ[K] \otimes_{\widetilde\QQ} \bar{U}(\PP^1)^-,
$$
where $K = K_0(\PP^1) \cong \ZZ^2$ is the $K$--group and the element
$C = K_{(0,1)}$ is central, see also \cite{HallEll}.
 Consider the group $K' := \ZZ \oplus \frac{1}{2} \ZZ$
and the algebra $\widetilde{DU}(\PP^1)$ obtained from
$DU(\PP^1)$ by adding two  central generators  $C^{\pm \frac{1}{2}} = K_{(0, \pm \frac{1}{2})}$ such that
$C^{\frac{1}{2}} C^{-\frac{1}{2}} = \mathbbm{1} = C^{-\frac{1}{2}} C^{\frac{1}{2}}$ and
$\bigl(C^{\frac{1}{2}})^2 = C$. For any $r \in \ZZ_{> 0}$ we set:
$
\widetilde{T}_r^{\pm}  = T_r^\pm  \cdot C^{\mp \frac{r}{2}} \quad
\mbox{and} \quad
\widetilde{\Theta}_r^{\pm}  = \Theta_r^\pm  \cdot C^{\mp \frac{r}{2}}.
$
Then we have:
$$
\Delta\bigl(\widetilde{T}_r^+\bigr) = \widetilde{T}_r^+  \otimes C^{-\frac{r}{2}}
+ C^{\frac{r}{2}} \otimes \widetilde{T}_r^+, \quad r \in \ZZ_{>0}
$$
and
$$
\Delta(L_n^+) = L_n^+  \otimes \mathbbm{1} + K C^n \otimes L_n^+
+ \sum\limits_{r = 1}^\infty \widetilde{\Theta}_r^+  K  C^{n - \frac{r}{2}} \otimes
L_{n-r}^+, \quad n \in \ZZ.
$$
Note that by Lemma \ref{L:someGreenprod} we have:
 $\bigl(\widetilde{\Theta}_r, \widetilde{T}_r\bigr) = \bigl(\Theta_r, T_r\bigr) =
\frac{\displaystyle [2r]}{\displaystyle r}$.
Using the relations of the Drinfeld double and rewriting the relations of
Theorem \ref{T:relationsinP1},  we obtain:
\begin{enumerate}\item $\bigl[K, \widetilde{T}_n^\pm\bigr] = 0 = \bigl[\widetilde{T}_n^\pm, \widetilde{T}_m^\pm\bigr]$ for all $m, n  \in \ZZ_{> 0}$;
\item $K L_n^\pm = v^{\mp 2} L_n^\pm K$ for all $n \in \ZZ$;
\item $\bigl[\widetilde{T}_r^\pm, L_n^\pm \bigr] = \frac{\displaystyle [2r]}{\displaystyle r} L^\pm_{n+r} C^{\mp \frac{r}{2}}$
       for all $n \in \ZZ$ and $r \in \ZZ_{> 0}$;
\item  $\bigl[L_n^\pm, \widetilde{T}_r^\mp\bigr] = \frac{\displaystyle [2r]}{\displaystyle r} L_{n-r}^\pm C^{\mp \frac{r}{2}}$
       for all $n \in \ZZ$ and $r \in \ZZ_{> 0}$;
\item $\bigl[\widetilde{T}_r^+, \widetilde{T}_s^-\bigr] = \delta_{r,s} \frac{\displaystyle
[2r]}{\displaystyle r}
\frac{\displaystyle C^{-r} - C^r}{\displaystyle v^{-1} - v}$, where $r, s \in \ZZ_{> 0}$;
\item Finally, we have:
$$
\bigl[L_n^+, L_m^-\bigr] =
\left\{
\begin{array}{ccc}
 \frac{\displaystyle v}{\displaystyle v- v^{-1}} \,\widetilde{\Theta}_{n-m}^+ K C^{\frac{m+n}{2}} & \mbox{if} & n > m, \\
0 & \mbox{if} & n =  m,\\
\frac{\displaystyle v}{\displaystyle v^{-1}-v} \,\widetilde{\Theta}_{m-n}^- K^{-1} C^{-\frac{m+n}{2}} & \mbox{if} & n < m.
\end{array}
\right.
$$
\end{enumerate}

\begin{definition}\label{D:DrinfeldNew}
Consider the $R$--algebra $U_v(\mathfrak{L sl}_2)$ generated by the elements
$X_n^\pm \, (n \in \ZZ)$, $H_r \, (r \in \ZZ \setminus \{0\})$, $C^{\pm \frac{1}{2}}$ and
$K^{\pm}$ subject to the following relations:
\begin{enumerate}
\item $C^{\frac{1}{2}}$ is central;
\item $K^\pm K^\mp = 1 = K^\mp K^\pm$, $C^{\frac{1}{2}} C^{-\frac{1}{2}} = 1 =
C^{-\frac{1}{2}} C^{\frac{1}{2}}$;
\item $\bigl[K, H_r\bigr] = 0$ for all $r \in \ZZ \setminus \{0\}$,  $K X_n^\pm  = v^{\mp 2} X_n^\pm  K$ for all $n \in \ZZ$;
\item We have Heisenberg-type relations
 $$
 \bigl[H_m, H_n\bigr] =
 \delta_{m+n, 0} \frac{[2m]}{m} \frac{C^{n} - C^{-n}}{v - v^{-1}}
 $$
  for all $m, n  \in  \ZZ \setminus \{0\}$;
\item We have Hecke-type relations
$$\bigl[H_r, X_n^\pm\bigr] = \pm \frac{\displaystyle [2r]}{\displaystyle r} X_{n+r}^\pm C^{\mp \frac{|r|}{2}}$$ for all
$n \in \ZZ$ and $r  \in  \ZZ \setminus \{0\}$;
\item $X_{m}^\pm  X_{n+1}^\pm  + X_{n}^\pm X_{m+1}^\pm =  v^{\pm 2} \bigl(X_{n+1}^\pm  X_{m}^\pm
 + X_{m+1}^\pm X_{n}^\pm)$ for all $m, n \in
\ZZ$;
\item
Finally, for all $m, n \in \ZZ$ we have:
$$
\bigl[X_m^+, X_n^-\bigr] = \frac{v}{v - v^{-1}}
\bigl(\Psi_{m+n}^+ C^{\frac{m-n}{2}} - \Psi_{m+n}^- C^{\frac{n-m}{2}}\bigr) K^{\mathsf{sign}(m+n)},
$$
where $\Psi^\pm_{\pm r} (r \ge 1)$  are given by the generating series
$$
1 + \sum\limits_{r=1}^\infty \Psi^\pm_{\pm r} t^r = \exp\bigl(\pm(v^{-1}-v)\sum\limits_{r = 1}^\infty H_{\pm r} t^{r}\bigr)
$$
and $\Psi^\pm_{\pm r} = 0$ for  $r < 0$.
\end{enumerate}
\end{definition}

\begin{remark}\label{R:formofrelatforloop}
Similarly to the case of $U_v(\widehat{\mathfrak{sl}}_2)$, our presentation of
$U_v(\mathfrak{L sl}_2)$ slightly differs from the conventional one, as used
in \cite{Lusztig, ChariPressley, JimboMiwa}. To pass to their notation,
one has to replace $v$ by $v^{-1}$ at the first step and then  replace
$v X_n^-$ by $- X_n^-$ for all $n \in \ZZ$ at the second, see also
Remark \ref{R:formofrelat}.
\end{remark}

\begin{proposition}
Let $U_q(\mathfrak{L sl}_2) = U_v(\mathfrak{L sl}_2) \otimes_{R} \widetilde{\QQ}_q$.
Then the map $\ev_q: U_q(\mathfrak{L sl}_2) \rightarrow  \widetilde{DU}(\PP^1)$  given by the rule:
$X_n^+ \mapsto L_n^+$, $X_n^- \mapsto L_{-n}^-$  for  $n \in \ZZ$,
$H_r \mapsto \widetilde{T}_r^+$ for $r \in \ZZ_{> 0}$ and
$H_r \mapsto -\widetilde{T}_{-r}^-$ for $r \in \ZZ_{< 0}$,
$\Psi_r^+ \mapsto \widetilde{\Theta}_r^+$ for $r \in \ZZ_{> 0}$ and
$\Psi_r^- \mapsto \widetilde{\Theta}_{-r}^-$ for $r \in \ZZ_{< 0}$, $K \mapsto K$ and
$C^{\frac{1}{2}} \mapsto C^{\frac{1}{2}}$, is an isomorphism of algebras.
\end{proposition}

\begin{proof} From the list of relations  of  Theorem \ref{T:relationsinP1} it follows that
the morphism $\ev_q$ is well-defined.
Next, consider the elements $$B_{\underline{m}',\, \underline{l}', \,  a,\,  b, \, \underline{m}'', \, \underline{l}''} =
\prod_{n \in \ZZ} [X_n^+]^{(m'_n)} \circ \prod_{r \in \ZZ_{>0}}  (H_r)^{l'_r} \circ K^a C^{\frac{b}{2}}
\circ
\prod_{n \in \ZZ} [X_n^-]^{(m''_n)} \circ \prod_{r \in \ZZ_{ < 0}}  (H_r)^{l''_r}
$$
of $U_v(\mathfrak{L sl}_2)$,
where $(m'_n)_{n \in \ZZ}, (m''_n)_{n \in \ZZ}, (l'_r)_{r \in \ZZ_{>0}},
(l''_r)_{r \in \ZZ_{>0}}$ run through the set of all sequences
of non-negative integers such that all but finitely many entries are zero, and
$a, b \in \ZZ$.
Using the defining relations
it is not difficult to show  that all the elements $B_{\underline{m}', \, \underline{l}', \, a, \,   b, \, \underline{m}'', \, \underline{l}''}$ generate $U_v(\mathfrak{L sl}_2)$ as $R$-module.
Indeed, observe first that any element of $U_v(\mathfrak{L sl}_2)$ can be written as an $R$--linear combination
of elements of the form $A_+ \cdot Z \cdot A_{-}$, where
$Z \in \bigl\langle C^{\pm \frac{1}{2}}, K^\pm \bigr\rangle$ and
$A_{\pm} \in \bigl\langle X_n^\pm, H_{\pm r} | \; n \in \mathbb{Z}, \, r \in \mathbb{Z}_{>0}
\bigr\rangle$. Next, using  Hecke-type relations of Definition \ref{D:DrinfeldNew},
we can write any element $A \in \bigl\langle X_n, H_{r},  C^{\pm \frac{1}{2}} \big| \;  n \in \mathbb{Z}, \, r \in \mathbb{Z}_{>0}
\bigr\rangle$ as an $R$--linear combination of elements of the form $V \cdot T$, where
$V \in \bigl\langle X_n | \; n \in \mathbb{Z}
\bigr\rangle$ and $T \in \bigl\langle H_{r},  C^{\pm \frac{1}{2}} \big| \;  r \in \mathbb{Z}_{>0}
\bigr\rangle$. Finally, it remains to note that for all $m, n \in \mathbb{Z}$ such that
$m > n$ we can write $X_m X_n$ as an $R$--linear combination of the
elements $X_n X_m, X_{n+1} X_{m-1}, \dots, X_{n+l} X_{m-l}$, where $l$ is the entier of
$\frac{m-n}{2}$.

Moreover, for any
$q \in \kP$ the elements
$\ev_q\bigl(B_{\underline{m}', \,  \underline{l}', \, a, \,  b, \,    \underline{m}'', \,   \underline{l}''}\bigr)$ are
linearly independent in $\widetilde{DU}(\PP^1)$. Hence, $ U_v(\mathfrak{L sl}_2)$ is free as an $R$--module
and $B_{\underline{m}', \,  \underline{l}', \,  a, \,   b, \,   \underline{m}'', \,   \underline{l}''}$  is its basis over $R$.
In particular, the morphism $\ev_q$ is an isomorphism for any $q \in \kP$.
\end{proof}

\begin{corollary} Let
$\widetilde{DU}_{\mathrm{gen}}(\PP^1):= \prod\limits_{q \in \kP} DU\bigl(\PP^1(\FF_q)\bigr)$
then the $R$--linear  map $$\mathsf{ev} = \prod\limits_{q \in \kP} {\mathsf{ev}_q}: U_v(\mathfrak{L sl}_2) \lar \widetilde{DU}_{\mathrm{gen}}(\PP^1)$$
is injective. Moreover, the elements $B_{\underline{m}', \,  \underline{l}', \,  a, \,   b, \,   \underline{m}'', \,  \underline{l}''}$ form a Poincar\'e-Birkhoff-Witt  (PBW)  basis of $U_v(\mathfrak{L sl}_2)$ viewed as an $R$--module.
\end{corollary}

\begin{remark}
The Hall algebra approach to a construction of PBW bases of quantum groups is due to Ringel
\cite{RingelPBW}. The case of $U_v^+(\mathfrak{L sl}_2)$ was considered by Baumann and Kassel
 \cite{BK}, see also \cite{Zhang}.  PBW bases
 of quantum loop algebras were
studied by Beck \cite{Beck} and Beck, Chari and Pressley \cite{BeckChariPressley}.
\end{remark}

\begin{remark}\label{R:actofSerreonP1}
Consider the functor $\AA = \kO_{\PP^1}(-2) \otimes \,-:  \hspace{1mm} D^b\bigl(\Coh(\PP^1)\bigr)
\rightarrow D^b\bigl(\Coh(\PP^1)\bigr).$
 Then it     induces  an automorphism of the algebra
$\widetilde{DU}(\PP^1)$ preserving the subalgebra $DU(\PP^1)$ and given by:
$$
L_n^{\pm} \stackrel{\AA}\lar L_{n \mp 2}^\pm, \;
T_r^\pm  \stackrel{\AA}\lar T_r^\pm,  \; \Theta_r^\pm  \stackrel{\AA}\lar \Theta_r^\pm, \;
\widetilde{T}_r^\pm  \stackrel{\AA}\lar \widetilde{T}_r^\pm, \;
\widetilde{\Theta}_r^\pm  \stackrel{\AA}\lar \widetilde{\Theta}_r^\pm, \;
  C^{\frac{1}{2}} \stackrel{\AA}\lar C^{\frac{1}{2}}, \;
K \stackrel{\AA}\lar K C^{-2}.
$$
In particular, it  corresponds to the algebra automorphism $\AA$
of  $U_v(\mathfrak{L sl}_2)$ such that
$$
X_n^{\pm} \stackrel{\AA}\lar X_{n \mp 2}^\pm, \;
H_r  \stackrel{\AA}\lar H_r, \; \Psi_r  \stackrel{\AA}\lar \Psi_r, \;
 C^{\frac{1}{2}} \stackrel{\AA}\lar C^{\frac{1}{2}}, \;
K \stackrel{\AA}\lar K C^{-2}.
$$
Moreover, the following diagram is commutative:
$$
\xymatrix
{U_v(\mathfrak{L sl}_2) \ar[rr]^{\AA} \ar[d]_{\ev} & & U_v(\mathfrak{L sl}_2) \ar[d]^{\ev} \\
\widetilde{DU}_{\mathrm{gen}}(\PP^1) \ar[rr]^{\AA} & & \widetilde{DU}_{\mathrm{gen}}(\PP^1).
}
$$

\end{remark}

\section{Categorification of Drinfeld-Beck isomorphism for $U_v\bigl(\widehat{\mathfrak{sl}}_2\bigr)$}

\noindent
Next, we elaborate a connection between the reduced Drinfeld doubles
$DC(\overr{Q})$ and $DU(\PP^1)$.

\begin{theorem}\label{T:Tilting}
Let $\kF = \kO_{\PP^1} \oplus \kO_{\PP^1}(1)$ and $B = \End_{\PP^1}(\kF)$.
Then  $$\FF := \RR \Hom_{\PP^1}\bigl(\kF,\,-\,\bigr): D^b\bigl(\Coh(\PP^1)\bigr)
\lar D^b\bigl(\mod-B\bigr)$$ is an equivalence of triangulated categories.
Identifying the category of right $B$--modules  with the category of representations of the Kronecker
quiver $\overr{Q}$ we have the following statements.
\begin{enumerate}
\item In the diagram of categories and functors
$$
\xymatrix{
D^b\bigl(\Coh(\PP^1)\bigr) \ar[rr]^{\FF} \ar[d]_{\kO_{\PP^1}( \mp 2) \otimes \,-\,}& & D^b\bigl(\Rep(\overr{Q})\bigr) \ar[d]^{(\RR\SS^\pm)^2}\\
D^b\bigl(\Coh(\PP^1)\bigr) \ar[rr]^{\FF} & & D^b\bigl(\Rep(\overr{Q})\bigr)
}
$$
both compositions are isomorphic.
\item $\FF\bigl(\kO_{\PP^1}(n)\bigr) \cong  P_{n}$ if $n\ge 0$ and $I_{-n-1}[-1]$ if $n < 0$.
\item $\FF$ induces an equivalence between the category $\Tor(\PP^1)$ of torsion coherent sheaves
on $\PP^1$ and the subcategory $\mathsf{Tub}(\overr{Q})$  of $\Rep(\overr{Q})$, which is the additive closure of the category of modules lying in the tubes.
\end{enumerate}
 \end{theorem}

 \begin{proof}
 The result  that $\FF$ is an equivalence of categories is
due to Beilinson \cite{Beilinson}, see also \cite{GL}. Next, the functor $\AA= \kO_{\PP^1}(-2) \otimes \,-\,$ is the Auslander-Reiten translate in  $D^b\bigl(\Coh(\PP^1)\bigr)$,
see \cite{BondalKapranov, ReitenBergh}. On the other hand,
by Theorem \ref{T:ReflFunct} the functor
$(\RR\SS^\pm)^2$
is the Auslander-Reiten translate in  $D^b\bigl(\Rep(\overr{Q})\bigr)$. Since $\FF$ is an equivalence, we have: $\FF \circ \AA \cong (\RR\SS^+)^2 \circ \FF$, see for example
\cite[Proposition I.2.3]{ReitenBergh}.

It is clear that the algebra  $B = \End_{\PP^1}(\kF)$
is \emph{isomorphic} to the path algebra of the Kronecker
quiver. However, in order to describe the images of the indecomposable objects of
$D^b\bigl(\Coh(\PP^1)\bigr)$ in a precise way, we have to specify this isomorphism.
Recall that  a choice of homogeneous coordinates $(x: y)$ on $\mathbb{P}^1$ fixes  two distinguished
sections  $x, y \in  \Hom_{\mathbb{P}^1}\bigl(\kO_{\PP^1}, \kO_{\PP^1}(1)\bigr)$ vanishing at $(0: 1)$ and $(1: 0)$ respectively.
Let $e_1$ and $e_2$  be the primitive idempotents of $A:= B^{\mathsf{op}}$
corresponding  to the identity endomorphisms of  $\kO_{\PP^1}$ and  $\kO_{\PP^1}(-1)$ respectively. Then we have: $\FF\bigl(\kO_{\PP^1}\bigr) \cong A e_1$ and
$\FF\bigl(\kO_{\PP^1}(-1)\bigr) \cong A e_2$.  We identify the path algebra of the Kronecker quiver
$\overr{Q}$ with the algebra $A$ in such a way that the sections $x$ and $y$ got identified
with the upper  and lower  arrows of $\overr{Q}$ under
$$
\Hom_{\mathbb{P}^1}\bigl(\kO_{\PP^1}, \kO_{\PP^1}(1)\bigr) \stackrel{\GG}\lar \Hom_A\bigl(A e_2, A e_1\bigr)
=   e_2 A e_1 \cong  e_2 \cdot  \kk \overr{Q} \cdot e_1.
$$
Note that  for all $n \in \mathbb{Z}$ the triangle
$$\kO_{\PP^1}(n-1) \xrightarrow{(-x \, y)}  \kO_{\PP^1}(n)^{\oplus 2}
\xrightarrow{
\left(\begin{smallmatrix} y \\ x\end{smallmatrix}\right)} \kO_{\PP^1}(n+1) \lar   \kO_{\PP^1}(n-1)[1]$$
is almost split in $D^b\bigl(\Coh(\PP^1)\bigr)$. Since   $\FF$ maps almost split triangles to almost split triangles,  Theorem \ref{T:indofKron}  implies   the formula
for the images  of  the lines bundles under the functor $\FF$.
Applying Theorem \ref{T:Dold}, Theorem \ref{T:indofKron} and  Theorem \ref{T:basicsonP1}, one can easily  deduce that $\FF$ restricts to
an equivalence of abelian categories  $\Tor(\PP^1)$ and $\mathsf{Tub}(\overr{Q})$.

This correspondence can be made more precise.
Let $\pi   = (a: b) \in \PP^1$ be a closed point of \emph{degree one} given by the homogeneous
form $p_{\pi}(x, y) =  a y - bx$.
 Then the unique simple   torsion sheaf $\kT_{1, \, \pi}$ supported at $\pi$ has
a locally free resolution
$$
0 \lar \kO_{\PP^1}(-1) \stackrel{p_\pi}\lar  \kO_{\PP^1} \lar \kT_{1, \, \pi} \lar 0.
$$
Since $\GG\bigl(\kO_{\PP^1}(-1)\bigr) =  P_0$ and  $\GG\bigl(\kO_{\PP^1}\bigr) = P_1$, we conclude
that $\GG\bigl(\kT_{1,\, \pi}\bigr) \cong T_{1, \, \pi}$. It implies that for any $n \in \mathbb{Z}_{>0}$
we have: $\GG\bigl(\kT_{n, \, \pi}\bigr)  \cong T_{n, \, \pi}$. The general case can be treated similarly.
 \end{proof}

\noindent
Applying Cramer's Theorem \ref{T:Cramer} we obtain:

 \begin{corollary}
 The assignment $DH(\PP^1) \stackrel{\FF}\to DH(\overr{Q})$ is an isomorphism of algebras.
 \end{corollary}

Our next goal is to show this isomorphism restricts on the isomorphism between
the reduced Drinfeld doubles of the composition subalgebras
$DU(\PP^1) \to DC(\overr{Q})$. For this it is convenient to consider a
functor $\GG: D^b\bigl(\Rep(\overr{Q})\bigr) \to D^b\bigl(\Coh(\PP^1)\bigr)$, which is
quasi-inverse to $\FF$.

\begin{theorem}\label{T:isoofDDoofcom}
The algebra isomorphism $\GG: DH(\overr{Q}) \to DH(\PP^1)$ restricts on the isomorphism
of the reduced Drinfeld doubles of the composition subalgebras
 $ \GG: DC(\overr{Q}) \rightarrow DU(\PP^1)$.
\end{theorem}

\begin{proof}
We use the notation $E_i = [S_i]^+$, $F_i = [S_i]^-$ and $K_i = K_{\bar{S}_i}$, $i = 1,2$
for the elements of the algebra $DC(\overr{Q})$. Then we have:
$$
E_1 \stackrel{\GG}\lar v^{-1}  L_{-1}^{-} K^{-1} C, \, E_2 \stackrel{\GG}\lar   L_0^+, \,
F_1 \stackrel{\GG}\lar v^{-1} L_{-1}^+ K C^{-1}, \, F_2 \stackrel{\GG}\lar   L_0^-,
K_1 \stackrel{\GG}\lar K^{-1} C,  K_2
\stackrel{\GG}\lar K.
$$
This implies that the image of the subalgebra $DC(\overr{Q})$ is contained in
$DU(\PP^1)$. Moreover, from the commutativity of the diagram
$$
\xymatrix
{DC(\overr{Q}) \ar[rr]^{\GG} \ar@{^{(}->}[d] & & DU(\PP^1) \ar@{_{(}->}[d]  \\
DH(\overr{Q}) \ar[rr]^{\GG} & & DH(\PP^1)
}
$$
we derive  that the algebra homomorphism $\GG:  DC(\overr{Q}) \to  DU(\PP^1)$ is injective.
To show it is also surjective, we have to prove that all the elements
$T_r^\pm$ and $\Theta_r^\pm$ are in the image of $\GG$ for all $r \ge 1$.
Note that for any pair of integers $n > m$ we have the following relation in  $DU(\PP^1)$:
$$
\bigl[L_n^+,  L_m^-\bigr] = \frac{v}{v - v^{-1}} \Theta_{n-m}^+ K C^m.
$$
This implies that $\Theta_{r}^\pm$ belong  to $\GG\bigl(DC(\overr{Q})\bigr)$ for all $r > 0$.
Hence, the elements $T_r^\pm$ belong to  $\GG\bigl(DC(\overr{Q})\bigr)$ for all $r > 0$. This shows the surjectivity of the map $\GG: DC(\overr{Q}) \to  DU(\PP^1)$.
\end{proof}

\noindent
As an application of the developed  technique, we get a shorter and (on our mind) more conceptual
 proof of the following formula,
obtained for the first time by Sz\'ant\'o \cite[Theorem 4.3]{Szanto}.

\begin{theorem}\label{T:Szantosformula}
In the Hall algebra of the Kronecker quiver $\overr{Q}$ we have for any $m, n
\in  \ZZ_{\ge 0}$:
$$
[I_m] \cdot  [P_n]  - v^2 [P_n] \cdot [I_m] =
\frac{v^{-(m+n+1)}}{v^{-1}-v} \sum\limits_{\substack{\pi_1, \dots, \pi_l \in \kQ; \, \pi_i \ne \pi_j  \,
1 \le i \ne j \le l \\ t_1, \dots, t_l \,:\,
\sum_{i=1}^l t_i \deg(\pi_i) =  m+n +1 }}
\prod_{i = 1}^l \bigl(1 - v^{2\deg(\pi_i)}\bigr)\bigl[T_{t_i, \pi_i}\bigr],
$$
where we sum over the set $\kQ$ of all  homogeneous prime ideals of height one  in the ring
$k[x, y]$.
In particular the left-hand side of this formula depends only on the sum $m+n$.
\end{theorem}

\begin{proof}
Let $m$ and $n$ be non-negative integers such that $m + n + 1 = r$.
Then we  have the following identity in the reduced Drinfeld double
$DU(\PP^1)$:
$$
L_{-m-1}^- L_n^+ - L_n^+ L_{-m-1}^- =
\frac{1}{q-1} \Theta_r^+ K C^{-m-1}.
$$
Note that the algebra homomorphism $\FF: DU(\PP^1) \rightarrow  DC(\overr{Q})$ acts as follows:
$$
L_{-m-1}^- \mapsto v [I_m]^+ K_1^{-m -1} K_2^{-m}, \quad
L_n^+ \mapsto [P_n]^+, \quad  K_{(1, -t)} = K C^{-t} \mapsto K_1^{-t} K_2^{-t +1},
$$
where $m, n \in \ZZ_{\ge 0}$.
It remains to observe that by Remark \ref{R:anotherformofTheta} we have:
$$
\overline{\Theta}_{r}:= \FF(\Theta_{r}) =  v^{-r} \sum\limits_{\substack{t_1, \dots, t_l \, : \,
\sum_{i=1}^l t_i \deg(\pi_i) =  r \\ \pi_1, \dots, \pi_l \in \kQ; \, \pi_i \ne \pi_j 1 \le i \ne j \le l}}
\prod_{i = 1}^l \bigl(1 - v^{2\deg(\pi_i)}\bigr)\bigl[T_{t_i, \pi_i}\bigr].
$$
In particular, we get the equality:
$$
v [I_m] K_1^{-m-1} K_2^{-m} [P_n] - v [P_n] [I_m] K_1^{-m-1} K_2^{-m} =
\frac{1}{q-1} \overline{\Theta}_r K_1^{-m-1} K_2^{-m}.
$$
Taking into account the fact that $K_1 K_2$ is central and
$K_1^{-1} [P_n] = v^{-2} [P_n] K_{1}^{-1}$, we end up precisely with Sz\'ant\'o's formula.
Our proof explains  the  conceptional meaning  of this equality: this formula  in the composition
subalgebra of $\Rep(\overr{Q})$ is a translation of a ``canonical'' relation in the reduced
Drinfeld double of $U(\PP^1)$.
\end{proof}

Another application of our approach is the following important result, which was stated by  Drinfeld
\cite{Drinfeld} and proven by
Beck \cite{Beck}, see also \cite{DingFrenkel} and \cite{Jing}.

\begin{theorem}\label{T:DrinfBeck}
We have an injective homomorphism of $R$--algebras
$\GG: U_v\bigl(\widehat{\mathfrak{sl}}_2\bigr) \lar U_v(\mathfrak{Lsl}_2)$ given by the following formulae:
$$
E_1 \stackrel{\GG}\lar v^{-1} X_{1}^- K^{-1} C, \,  E_2 \stackrel{\GG}\lar X_0^+, \,
F_1 \stackrel{\GG}\lar v^{-1} X_{-1}^+ K C^{-1}, \, F_2 \stackrel{\GG}\lar X_0^-,
K_1 \stackrel{\GG}\lar K^{-1} C,   K_2 \stackrel{\GG}\lar K.
$$
Its image is the subalgebra of $U_v(\mathfrak{Lsl}_2)$ generated by the elements
$X_n^\pm, H_r, C^\pm$ and $K^\pm$. Moreover, the following diagram is commutative for any
$m \in \ZZ$:
$$
\xymatrix{
U_v\bigl(\widehat{\mathfrak{sl}}_2\bigr) \ar[rr]^{\AA^m} \ar[d]_{\GG} & &
U_v\bigl(\widehat{\mathfrak{sl}}_2\bigr) \ar[d]^{\GG}\\
U_v(\mathfrak{Lsl}_2) \ar[rr]^{\AA^m} & & U_v(\mathfrak{Lsl}_2).
}
$$
\end{theorem}

\begin{proof}
This result follows directly from the commutativity of the following diagram:
$$
\xymatrix{
U_v\bigl(\widehat{\mathfrak{sl}}_2\bigr) \ar[rrr]^{\AA^\pm} \ar[ddd]_{\GG}  \ar[dr]_\ev & & &
U_v\bigl(\widehat{\mathfrak{sl}}_2\bigr) \ar[ddd]^{\GG} \ar[dl]^\ev\\
& {DC}_{\mathrm{gen}}(\overr{Q}) \ar[r]^{\AA^\pm} \ar[d]_{\GG} &  {DC}_{\mathrm{gen}}(\overr{Q}) \ar[d]^{\GG} &
\\
& \widetilde{DU}_{\mathrm{gen}}(\PP^1) \ar[r]^{\AA^\pm} &  \widetilde{DU}_{\mathrm{gen}}(\PP^1) &
\\
U_v(\mathfrak{Lsl}_2) \ar[rrr]^{\AA^\pm} \ar[ur]^\ev & & & U_v(\mathfrak{Lsl}_2), \ar[ul]_\ev
}
$$
which is obtained by  patching together  the diagrams, constructed in
 Corollary \ref{C:ReflonCompAlg},
Remark \ref{R:actofSerreonP1} and Theorem \ref{T:isoofDDoofcom}.
\end{proof}

\begin{remark}
The algebra homomorphism  $\GG: U_v\bigl(\widehat{\mathfrak{sl}}_2\bigr) \to
U_v(\mathfrak{Lsl}_2)$ is \emph{not surjective} because the elements $C^{\pm \frac{1}{2}}$
 do not belong to the image of $\GG$.
\end{remark}

\begin{remark}
In order to pass to the ``conventional form'' of the Drinfeld-Beck (iso)morphism, one has
to apply the identifications of  $U_v(\widehat{\mathfrak{sl}}_2)$ and $U_v(\mathfrak{L sl}_2)$ described
in Remark \ref{R:formofrelat} and  Remark \ref{R:formofrelatforloop}. In that terms, the ``categorical
isomorphism'' $\GG$  of Theorem \ref{T:DrinfBeck}
is equal to the composition  of  the conventional one (as can be found for instance in
\cite{JimboMiwa})  combined with the  automorphism of $U_v(\widehat{\mathfrak{sl}}_2)$ given by the
rule $E_1 \mapsto - C^{-1}E_1$, $F_1 \mapsto -CF_1$ and leaving the remaining generators unchanged.
\end{remark}

\medskip

\noindent
As a final application of the technique of Hall algebras to the study of the quantum affine algebra
$U_v\bigl(\widehat{\mathfrak{sl}}_2\bigr)$, we shall reprove  several known results on its
 integral form.

\begin{definition}
The integral form $U_v^{\mathrm{int}}\bigl(\widehat{\mathfrak{sl}}_2\bigr)$ of the
quantum affine algebra $U_v\bigl(\widehat{\mathfrak{sl}}_2\bigr)$
is the $\QQ[v, v^{-1}]$ subalgebra of $U_v\bigl(\widehat{\mathfrak{sl}}_2\bigr)$
generated by  $E_i^{(n)} = \frac{\displaystyle E_i^n}{\displaystyle [n]!}$,
$F_i^{(n)} = \frac{\displaystyle F_i^n}{\displaystyle [n]!}$ for all  $n \in \ZZ_{\ge 0}$, $i = 1,2$,
$K_1$ and $K_2$.
\end{definition}

\noindent
The following result is well-known, see for example \cite{Lusztig}.

\begin{theorem}
The algebra $U_v^{\mathrm{int}}\bigl(\widehat{\mathfrak{sl}}_2\bigr)$ is a Hopf algebra over
$\QQ[v, v^{-1}]$ and we
have: $$U_v^{\mathrm{int}}\bigl(\widehat{\mathfrak{sl}}_2\bigr) \otimes_{\QQ[v, v^{-1}]} R =
U_v\bigl(\widehat{\mathfrak{sl}}_2\bigr).$$
\end{theorem}

\noindent
Let $U_v^{\mathrm{int}}\bigl({\mathfrak{L sl}}_2\bigr) := \GG\bigl(U_v^{\mathrm{int}}\bigl(\widehat{\mathfrak{sl}}_2\bigr)\bigr)$. As an application
of our approach, we show that certain elements of
 $U_v\bigl({\mathfrak{L sl}}_2\bigr)$ actually belong to $U_v^{\mathrm{int}}\bigl({\mathfrak{L sl}}_2\bigr)$. First note the following well-known fact, see e.g.~\cite{IntForm}.

\begin{lemma}
The elements $X_n^{\pm(m)} \in U_v\bigl({\mathfrak{L sl}}_2\bigr)$  belong  to the integral form
$U_v^{\mathrm{int}}\bigl({\mathfrak{L sl}}_2\bigr)$ for all $n \in \ZZ$ and $m \in \ZZ_{>0}$.
\end{lemma}

\begin{proof}
Let $\cA = \Rep(\overr{Q})$ and $X \in \mathrm{Ob}(\cA)$ be an object such that
$\End_{\cA}(X) = \kk$ and $\Ext^1_{\cA}(X, X) = 0$. Then for any $n \in \ZZ_{\ge 0}$
we have the following equality  in the Hall algebra $H(\cA)$:
$$
\bigl[X^{\oplus n}\bigr] = v^{n(n-1)} \frac{[X]^n}{[n]!} = v^{n(n-1)} [X]^{(n)}.
$$
In particular, for any $i \in \{1,2\}$,
$n \in \mathbb{Z}_{> 0}$ and and $q \in \kP$ we have:
$\ev_q\bigl(E_i^{(n)}\bigr) = v^{n(1-n)} \bigl[S_i^{\oplus n}\bigr]$. Our next aim is to show
 the automorphisms $\SS^\pm$ of the algebra $U_v\bigl(\widehat{\mathfrak{sl}}_2\bigr)$
preserve the subalgebra $U_v^{\mathrm{int}}\bigl(\widehat{\mathfrak{sl}}_2\bigr)$. By  Corollary \ref{C:ReflonCompAlg}
it is sufficient to check that for any $n \in \ZZ_{>0}$  we have:
$\ev_q^{-1}\bigl([P_1^{\oplus n}]\bigr)$ and $\ev_q^{-1}\bigl([I_1^{\oplus n}]\bigr)$ belong
to the subalgebra $U_v^{\mathrm{int}, +}\bigl(\widehat{\mathfrak{sl}}_2\bigr)$. To show this, we
use the following trick.
For a  representation $X  = \xymatrix{
U  \ar@/^/[rr]^{A} \ar@/_/[rr]_{B} & & V}
$
of the Kronecker quiver $\overr{Q}$ we denote  $r(X) := \dim_{\kk}\bigl(\Im(A) + \Im(B)\bigr)$.
For any $n \in \mathbb{Z}_{>0}$ and $0 \le r \le 2n$ consider
$$
\widetilde{\mathbbm{1}}_{(n, 2n)}^{(r)} \; := \;
\sum\limits_{ \substack{[X] \in J: \; \underline{\dim}(X) = (n, 2n) \\
 r(X) = r}}
 [X] \; \; \in \; \; H(\cA).
$$
Then  $[P_1^{\oplus n}] = \widetilde{\mathbbm{1}}_{(n, 2n)}^{(2n)}$ and  for all
$a, b \in \ZZ_{\ge 0}$ such that $a + b = 2n$ we have the following identity:
$$
[S_2^{\oplus a}] \circ [S_1^{\oplus n}] \circ [S_2^{\oplus b}] =
v^{-ab + 2nb} \sum\limits_{r = 0}^b \Bigl|\mathsf{Gr}_\kk\bigl(b-r, 2n-r\bigr)\Bigr| \, \widetilde{\mathbbm{1}}_{(n, 2n)}^{(r)},
$$
where $\Bigl|\mathsf{Gr}_\kk\bigl(b-r, 2n-r\bigr)\Bigr| = v^{(b-2n)(b-r)}\frac{\displaystyle [2n-r]!}{\displaystyle [b-r]! [2n-b]!}$ is the
number of points  of  the Grassmanian $\mathsf{Gr}_\kk\bigl(b-r, 2n-r\bigr)$. From this formula one can deduce
that
$$
[P_1]^{(n)} = \sum\limits_{a+b=2n} (-1)^a v^{-b} [S_2]^{(a)} \circ [S_1]^{(n)} \circ [S_2]^{(b)}.
$$
In a similar way, one can show that
$$
[I_1]^{(n)} = \sum\limits_{a+b=2n} (-1)^a v^{-b} [S_1]^{(a)} \circ [S_2]^{(n)} \circ [S_1]^{(b)}.
$$
From the invariance of $U_v^{\mathrm{int}}\bigl(\widehat{\mathfrak{sl}}_2\bigr)$ under the
action of $\SS^\pm$  it
 also follows  that for any indecomposable pre-projective or pre-injective
object $X \in \mathrm{Ob}(\cA)$ and $n \in \ZZ_{>0}$,  the element
$[X^{\oplus n}]$ belongs to $C_{\mathrm{gen}}(\overr{Q})$
and lies  in the image of the homomorphism
$\ev: U_v^{\mathrm{int},+}(\widehat{\mathfrak{sl}}_2) \rightarrow C_{\mathrm{gen}}(\overr{Q})$.
Since for any $c \in \ZZ$ and $n \in \ZZ_{\ge 0}$ the vector bundle
$\kO_{\PP^1}(c)^{\oplus n}$ is isomorphic to $\GG(X^{\oplus n})$ for an appropriate
shift of an indecomposable  pre-projective or pre-injective module $X$, Theorem \ref{T:DrinfBeck} yields  the claim.
\end{proof}

\begin{lemma}\label{L:someidforcharf}
For any pair of  non-negative integers $(a ,b)$,  the  element
$$
\mathbbm{1}_{(a, b)} = \sum\limits_{[X] \in \cI \, : \, \bar{X} = (a, b)} [X] \in H(\overr{Q})
$$
belongs to the composition subalgebra $C(\overr{Q})$.
\end{lemma}

\begin{proof}
It  follows from the equality
$
\mathbbm{1}_{(a, b)} =  v^{-2ab} [S_1^{\oplus a}] \circ  [S_2^{\oplus b}].
$
\end{proof}

\begin{corollary}
For any pair of  non-negative integers $(a ,b)$  we have a well-defined element $\mathbbm{1}_{(a, b)}\in C_{\mathrm{gen}}(\overr{Q})$
belonging to the image of
$\ev: U_v^{\mathrm{int},+}(\widehat{\mathfrak{sl}}_2) \rightarrow C_{\mathrm{gen}}(\overr{Q}).$
\end{corollary}

\begin{lemma}\label{L:semistinHall}
Let $\cA = \Rep(\overr{Q})$,
$\HH = \{ a \in \CC \,| \, \mathrm{Im}(a) >0\}$ and $Z: K_0(\cA)
\to \CC$ be any additive group homomorphism such that for any
non-zero object $X$ of $\cA$ we have: $Z(\bar{X}) \in \HH$. For any $\alpha \in K_0(\cA)$ denote
$$
\mathbbm{1}^{\mathrm{ss}}_\alpha  = \mathbbm{1}^{\mathrm{ss}}_{\alpha, Z} := \sum\limits_{[X] \in \cI: \, X \in \cA^{\mathrm{ss}}_\alpha}
[X],
$$
where $\cA^{\mathrm{ss}}_\alpha$ is the category of semi-stable objects of class $\alpha$ with
respect to the stability condition $Z$, see \cite{Rudakov} for the definition.
  Then we have:
$\mathbbm{1}^{\mathrm{ss}}_\alpha \in C(\overr{Q})$.
\end{lemma}

\begin{proof}
First note that the existence and uniqueness of the Harder-Narasimhan filtration  \cite{Rudakov}
of an object of our abelian category $\cA$ implies the following
identity for an arbitrary class $\alpha \in K_0(\cA)$ and a given stability condition $Z$:
$$
\mathbbm{1}_\alpha = \mathbbm{1}^{\mathrm{ss}}_\alpha +
\sum\limits_{t\ge  2}  \sum\limits_{\substack{\alpha_1 + \dots + \alpha_t = \alpha \\
\mu(\alpha_1) \ge \dots \ge \mu(\alpha_t) }}
v^{\sum_{i < j} \langle \alpha_i \alpha_j \rangle}
\mathbbm{1}^{\mathrm{ss}}_{\alpha_1} \circ \dots \circ \mathbbm{1}^{\mathrm{ss}}_{\alpha_t}.
$$
Since the expression on the right-hand side is a finite sum, by induction
we obtain that for all classes $\alpha \in K_0(\cA)$ the element
$\mathbbm{1}^{\mathrm{ss}}_\alpha$ belongs to the subalgebra
of $H(\overr{Q})$ generated by all the elements $\bigl\{\mathbbm{1}_\beta\bigr\}_{\beta \in K_0(\cA)}$.
According to Lemma \ref{L:someidforcharf},   this algebra coincides with the composition
subalgebra $C(\overr{Q})$, what implies the claim.
\end{proof}

\begin{remark}\label{R:ReinekeInv}
A result  of Reineke \cite[Theorem 5.1]{Reineke} provides an explicit formula expressing
the elements $\bigl\{\mathbbm{1}^{\mathrm{ss}}_\alpha\bigr\}$ via $\bigl\{\mathbbm{1}_\beta\bigr\}$ for an arbitrary
stability function $Z$:
$$
\mathbbm{1}^{\mathrm{ss}}_{\alpha, Z} = \mathbbm{1}_\alpha +
\sum\limits_{t\ge  2} (-1)^{t-1}  \sum\limits_{\substack{\alpha_1 + \dots + \alpha_t = \alpha: \,
\forall \,  1 \le s \le t-1 \\
\mu(\alpha_1 + \dots + \alpha_s) > \mu(\alpha)}}
v^{\sum_{i < j} \langle \alpha_i \alpha_j \rangle}
\mathbbm{1}_{\alpha_1} \circ \dots \circ \mathbbm{1}_{\alpha_t}.
$$
Hence, the  element
$\mathbbm{1}^{\mathrm{ss}}_\alpha  = \mathbbm{1}^{\mathrm{ss}}_{\alpha, Z} \in C_{\mathrm{gen}}(\overr{Q})$ belongs to the
image of the algebra homomorphism $\ev: U_v^{\mathrm{int},+}(\widehat{\mathfrak{sl}}_2) \to  C_{\mathrm{gen}}(\overr{Q})$ for any class $\alpha \in K_0(\cA)$ and a given   stability condition
$Z: K_0(\cA) \to \HH$.
\end{remark}

\begin{lemma}\label{L:semistinCompAlg}
For any $r \in \ZZ_{>0}$ the following element of the Hall algebra $H(\overr{Q})$
$$
\widetilde{\mathbbm{1}}_{(r, r)} =  \sum\limits_{X \in \mathsf{Tub}(\overr{Q})\,:\, \bar{X} = (r, r)} [X]
$$
belongs to the composition algebra $C(\overr{Q})$. Moreover,
it determines  an element
of $C_{\mathrm{gen}}(\overr{Q})$ belonging to the
image of the homomorphism  $\ev: U_v^{\mathrm{int},+}(\widehat{\mathfrak{sl}}_2) \rightarrow   C_{\mathrm{gen}}(\overr{Q})$.
\end{lemma}

\begin{proof}
Consider the stability condition on the category $\Rep(\overr{Q})$ defined  by the
function $Z: K_0\bigl(\Rep(\overr{Q})\bigr) \to \RR^2$ given by the
rule $Z(m \bar{S}_1 + n \bar{S}_2) = (m-n, m+n)$.
Then the class  of a  pre-projective objects has  the form $(-1, l)$, $l \in \ZZ_{>0}$,
the class  of a  pre-injective representations has the form  $(1, l)$, $l \in \ZZ_{>0}$,
whereas the classes of the tubes have the form $(0, l)$, $l \in \ZZ_{>0}$.
Recall that
\begin{itemize}
\item an object of an abelian category  is semi-stable if and only if all its direct summands are semi-stable  with the same slope;
\item  any non-semi-stable object
can be destabilized by an indecomposable one;
\item there are no morphisms from a pre-injective object of $\Rep(\overr{Q})$
to an object from a tube.
\end{itemize}
Hence, $\widetilde{\mathbbm{1}}_{(r,r)} = \mathbbm{1}^{\mathrm{ss}}_{(r, r)}$
for all $r \in \ZZ_{>0}$. Applying Lemma \ref{L:semistinHall} and
Remark \ref{R:ReinekeInv}, we get the claim.
\end{proof}

\begin{remark}
Consider the standard stability condition on the category $\Coh(\PP^1)$ given
by the function $Z = (\rk, -\deg)$. Then it determines a stability condition
on the derived category $D^b\bigl(\Coh(\PP^1)\bigr)$ in the sense of Bridgeland \cite{Bridgeland}
 such that any indecomposable object of $D^b\bigl(\Coh(\PP^1)\bigr)$ is $Z$-semi-stable. In particular,
all  objects of the category $\Tor(\PP^1)$ are semi-stable
of slope $0$.
The stability condition on $\Rep(\overr{Q})$ used in the proof of Lemma \ref{L:semistinCompAlg}
induces  a stability condition on the derived category $D^b\bigl(\Rep(\overr{Q})\bigr)$.
In the notations of \cite{Bridgeland}, this stability condition is $\widetilde{\GL}(2, \RR)^+$-equivalent
to the standard stability condition on $D^b\bigl(\Coh(\PP^1)\bigr)$.

 The technique of stability conditions  also plays a key role
in our subsequent paper \cite{BSWeightLine} on the composition Hall algebra of a weighted projective line.
\end{remark}

We conclude this section with a new proof of the following proposition, which was obtained for the first time by Chari and Pressley in \cite{IntForm}.

\begin{proposition}
 Consider the family
 of  elements $\bigl\{P_r\bigr\}_{r \ge 1}$  of $U_v\bigl({\mathfrak{L sl}}_2\bigr)$ defined by the following generating series:
$$
1 + \sum\limits_{r = 1}^\infty P_r C^{-\frac{r}{2}} t^r =
\exp\Bigl(\sum\limits_{r = 1}^\infty \frac{\Psi_r}{[r]} t^r\Bigr).
$$
Then $P_r$ belong to the algebra
$U_v^{\mathrm{int}}\bigl({\mathfrak{L sl}}_2\bigr)$ for all $r \in \ZZ_{> 0}$.
\end{proposition}

\begin{proof}
The elements $P_r$ have a clear meaning in the language of the Hall algebra $H(\PP^1)$. Indeed,
by Definition \ref{D:compalggener} we have:
$\ev_q(P_r) =
\mathbbm{1}_{r\delta}
$ for any $q \in \kP$.

On the other side, Theorem \ref{T:Tilting} implies that
$\mathbbm{1}_{r\delta} = \GG\bigl(\widetilde{\mathbbm{1}}_{(r,r)}\bigr)$ for all $r \in \ZZ_{> 0}$.
By Lemma \ref{L:semistinCompAlg} we know that $\widetilde{\mathbbm{1}}_{(r,r)}$
belongs to the image of the algebra homomorphism $\ev: U_v^{\mathrm{int},+}(\widehat{\mathfrak{sl}}_2)
\rightarrow DC_{\mathrm{gen}}(\overr{Q})$. Hence,
$\mathbbm{1}_{r\delta}$ belongs to the image of the algebra homomorphism
$\ev: U_v^{\mathrm{int}}({\mathfrak{L sl}}_2) \rightarrow
DU_{\mathrm{gen}}(\PP^1)$. By Theorem \ref{T:DrinfBeck},  the element $P_r$ belongs
to the algebra $U_v^{\mathrm{int}}({\mathfrak{L sl}}_2)$ for all $r \in \ZZ_{> 0}$, too.
\end{proof}

\end{document}